# FUNCTIONAL LARGE DEVIATIONS FOR MULTIVARIATE REGULARLY VARYING RANDOM WALKS

By Henrik Hult,[1] Filip Lindskog,[2] Thomas Mikosch[3]
and Gennady Samorodnitsky[4]

*Cornell University, KTH, University of Copenhagen and Cornell University*

*Dedicated to the memory of Alexander V. Nagaev*

We extend classical results by A. V. Nagaev [*Izv. Akad. Nauk UzSSR Ser. Fiz.–Mat. Nauk* **6** (1969) 17–22, *Theory Probab. Appl.* **14** (1969) 51–64, 193–208] on large deviations for sums of i.i.d. regularly varying random variables to partial sum processes of i.i.d. regularly varying vectors. The results are stated in terms of a heavy-tailed large deviation principle on the space of càdlàg functions. We illustrate how these results can be applied to functionals of the partial sum process, including ruin probabilities for multivariate random walks and long strange segments. These results make precise the idea of *heavy-tailed large deviation heuristics*: in an asymptotic sense, only the largest step contributes to the extremal behavior of a multivariate random walk.

**1. Introduction and background.** The notion of regular variation is fundamental in various fields of applied probability. It serves as domain of attraction condition for partial sums of i.i.d. random vectors [26] or for component-wise maxima of vectors of i.i.d. random vectors [25], and it occurs in a natural way for the finite-dimensional distributions of the stationary solution to stochastic recurrence equations (see [11, 15]), including ARCH

Received November 2004; revised May 2005.
[1]Supported by the Swedish Research Council.
[2]Supported in part by Credit Suisse, Swiss Re and UBS through RiskLab, Switzerland.
[3]Supported in part by MaPhySto, the Danish research network for mathematical physics and stochastics, and the Danish Research Council (SNF) Grant 21-01-0546.
[4]Supported in part by NSF Grant DMS-03-03493 and NSA Grant MSPF-02G-183 at Cornell University.

and GARCH processes; see [2] and Section 8.4 in [8]. To start with, we consider an $\mathbb{R}^d$-valued vector **X**. We call it *regularly varying* if there exists a sequence $(a_n)$ of positive numbers such that $a_n \uparrow \infty$ and a nonnull Radon measure $\mu$ on the $\sigma$-field $\mathcal{B}(\overline{\mathbb{R}}^d \setminus \{\mathbf{0}\})$ of the Borel sets of $\overline{\mathbb{R}}^d \setminus \{\mathbf{0}\}$ such that $\mu(\overline{\mathbb{R}}^d \setminus \mathbb{R}^d) = 0$ and

(1.1) $$n\mathrm{P}(a_n^{-1}\mathbf{X} \in \cdot) \xrightarrow{v} \mu(\cdot),$$

where $\xrightarrow{v}$ denotes vague convergence on $\mathcal{B}(\overline{\mathbb{R}}^d \setminus \{\mathbf{0}\})$. We refer to [14] and [24, 25] for the concept of vague convergence. It can be shown that the above conditions on the distribution of **X** necessarily imply that $\mu(tA) = t^{-\alpha}\mu(A)$ for some $\alpha > 0$, all $t > 0$ and any Borel set $A$. Therefore, we also refer to *regular variation with index $\alpha$* in this context.

Definition (1.1) of regular variation has the advantage that it can be extended to random elements **X** with values in a separable Banach space (e.g., [1]) or certain linear metric spaces. Recently, de Haan and Lin [12] have used regular variation of stochastic processes with values in the space of continuous functions on $[0,1]$ to prove weak convergence results for the extremes of regularly varying processes with continuous sample paths. They also considered regular variation for stochastic processes with values in the Skorokhod space $\mathbb{D} = \mathbb{D}([0,1], \mathbb{R}^d)$ of $\mathbb{R}^d$-valued càdlàg functions on $[0,1]$, equipped with the $J_1$-topology (see [3]) very much in the same way as (1.1). This idea was taken up by Hult and Lindskog [13]. They characterized regular variation of càdlàg processes by regular variation of their finite-dimensional distributions in the sense of (1.1) and a relative compactness condition in the spirit of weak convergence of stochastic processes; see [3]. Then, not surprisingly, one can derive a continuous mapping theorem for regularly varying stochastic processes and apply it to various interesting functionals, including suprema of Lévy and Markov processes with weakly dependent increments.

In this paper we continue the investigations started by Hult and Lindskog [13] in a different direction. As a matter of fact, the notion of regular variation as defined in (1.1) is closely related to large deviation results for processes with heavy-tailed margins. Such results have been proved since the end of the 1960s by, among others, A. V. Nagaev [19, 20], S. V. Nagaev [21] and Cline and Hsing [5] for various one-dimensional settings; see Section 8.6 in [8] and [18] for surveys on the topic. In the mentioned papers it was shown for a random walk $S_n = Z_1 + \cdots + Z_n$ of i.i.d. random variables $Z_i$ that relations of the type

(1.2) $$\sup_{x \geq \lambda_n} \left| \frac{\mathrm{P}(S_n > x)}{n\mathrm{P}(Z_1 > x)} - 1 \right| \to 0$$

hold for suitable sequences $\lambda_n \to \infty$ and heavy-tailed distributions of $Z_i$. For example, S. V. Nagaev [21] showed that (1.2) holds for i.i.d. centered random variables $Z_i$ which are regularly varying with index $\alpha > 2$, where the sequence $(\lambda_n)$ can be chosen as $\lambda_n = a\sqrt{n \log n}$ for any $a > \sqrt{\alpha - 2}$. As a matter of fact, results of type (1.2) also hold for $Z_i$'s with a subexponential distribution. The latter class of distributions is wider than the class of regularly varying distributions. For our purposes, we will focus on regularly varying $Z_i$'s with index $\alpha > 0$. Then it follows from (1.2), using the uniform convergence theorem for regularly varying functions (see [4]), that

$$\sup_{x \geq 1} \left| \frac{\mathrm{P}(\lambda_n^{-1} S_n \in (x, \infty))}{n \mathrm{P}(Z_1 > \lambda_n)} - x^{-\alpha} \right| \to 0.$$

Motivated by this, we say that the partial sum process $\mathbf{S}_n = \mathbf{Z}_1 + \cdots + \mathbf{Z}_n$ of i.i.d. $\mathbb{R}^d$-valued regularly varying random vectors $\mathbf{Z}_i$ satisfies a *large deviation principle* if there exist sequences $\gamma_n, \lambda_n \uparrow \infty$ and a nonnull Radon measure $\mu$ on $\mathcal{B}(\overline{\mathbb{R}}^d \setminus \{\mathbf{0}\})$ such that

$$(1.3) \qquad \gamma_n \mathrm{P}(\lambda_n^{-1} \mathbf{S}_n \in \cdot) \xrightarrow{v} \mu(\cdot).$$

Similarly to the notion of regular variation, the latter definition allows one to extend large deviation principles from $\mathbb{R}^d$-valued sequences $(\mathbf{S}_n)$ to sequences of stochastic processes $(\mathbf{X}^n)$ with values in $\mathbb{D}$. This extension can be handled in the same way as for regular variation: one can give a criterion for a large deviation principle in terms of large deviation principles for the finite-dimensional distributions of the sequence $(\mathbf{X}^n)$ in combination with a relative compactness condition. As a consequence, one can derive a continuous mapping theorem.

The hard part of the proofs is to show the large deviation principle for the sequence $(\mathbf{X}^n)$. However, for the partial sums $\mathbf{S}_n$ of i.i.d. regularly varying $\mathbb{R}^d$-valued $\mathbf{Z}_i$'s, this is a relatively straightforward task. We show in Theorem 2.1 that a functional analogue to (1.3) with limiting measure $m$ holds for the $\mathbb{D}$-valued suitably centered processes $(\mathbf{S}_{[nt]})_{t \in [0,1]}$ with $\gamma_n = [n\mathrm{P}(|\mathbf{Z}| > \lambda_n)]^{-1}$. If the index of regular variation $\alpha > 1$, we may choose $\lambda_n = n$. The limiting measure $m$ is concentrated on step functions with one step. The interpretation is that, for large $n$, the process $\lambda_n^{-1} \mathbf{S}_{[n \cdot]}$ behaves like a step function with one step. As a consequence, we determine, in Theorem 3.1, the asymptotic behavior of the probability

$$\psi_u(A) = \mathrm{P}(\mathbf{S}_n - \mathbf{c}n \in uA \text{ for some } n \geq 1)$$

as $u \to \infty$. Here the steps $\mathbf{Z}_i$ are regularly varying with index $\alpha > 1$ and $\mathrm{E}(\mathbf{Z}_i) = \mathbf{0}$. Moreover, $\mathbf{c} \neq \mathbf{0}$ is a vector and $A$ is a set bounded away from some narrow cone in the direction $-\mathbf{c}$. The probability $\psi_u(A)$ may be interpreted as a multivariate ruin probability; ruin occurs when the random walk with drift $-\mathbf{c}$ hits the set $A$. If $\mu$ denotes the limiting measure in (1.3) of the random walk, then

$$\mu^*(A^\circ) \leq \liminf_{u\to\infty} \frac{\psi_u(A)}{u\mathrm{P}(|\mathbf{Z}|>u)}$$
$$\leq \limsup_{u\to\infty} \frac{\psi_u(A)}{u\mathrm{P}(|\mathbf{Z}|>u)} \leq \mu^*(\overline{A}),$$

where $A^\circ$ and $\overline{A}$ are the interior and closure of $A$, respectively, and for any set $B$,

$$\mu^*(B) = \int_0^\infty \mu(\mathbf{c}v + B_{\mathbf{c}})\,dv, \qquad B_{\mathbf{c}} = \{\mathbf{x} + \mathbf{c}t, \mathbf{x} \in B, t \geq 0\}.$$

For more details, see Section 3.

The functional large deviation result also applies to the asymptotic behavior of long strange segments of a random walk (see Section 4). Suppose $\alpha > 1$ and $\mathrm{E}(\mathbf{Z}_i) = \mathbf{0}$. For a set $A \in \mathcal{B}(\mathbb{R}^d)$ bounded away from $\mathbf{0}$, let

$$R_n(A) = \sup\{k : \mathbf{S}_{i+k} - \mathbf{S}_i \in kA \text{ for some } i \in \{0, \ldots, n-k\}\}.$$

A segment of length $R_n(A)$ is called a long strange segment. The name is motivated by observing that $R_n(A)$ is the length of an interval over which the sample mean is "far away" from the true mean. We show, in Theorem 4.1, that, for every $t \in (0,1)$ and $A \in \mathcal{B}(\mathbb{R}^d)$ bounded away from $\mathbf{0}$,

$$\mu(A^\circ(t)) \leq \liminf_{n\to\infty} \frac{\mathrm{P}(n^{-1}R_n(A) > t)}{n\mathrm{P}(|\mathbf{Z}|>n)}$$
$$\leq \limsup_{n\to\infty} \frac{\mathrm{P}(n^{-1}R_n(A) > t)}{n\mathrm{P}(|\mathbf{Z}|>n)} \leq \mu(A^*(t)),$$

where

$$A^*(t) = \bigcup_{t \leq s \leq 1} s\overline{A}, \qquad A^\circ(t) = \bigcup_{t < s \leq 1} sA^\circ.$$

In particular, if $A$ is an increasing set (i.e., $t\mathbf{x} \in A$ for $\mathbf{x} \in A$, $t \geq 1$) with $\mu(\partial A) = 0$, this simplifies to

$$\lim_{n\to\infty} \frac{\mathrm{P}(n^{-1}R_n(A) > t)}{n\mathrm{P}(|\mathbf{Z}|>n)} = t^{-\alpha}\mu(A).$$

From this result we derive, in Theorem 4.2, the weak limit of $(a_n^{-1}R_n(A))$, where $(a_n)$ is the sequence associated with the regularly varying $\mathbf{Z}_i$'s in (1.1).

We want to mention that some of the technical issues encountered in the proofs in this paper arise when switching from the discrete time random walk to the continuous time limit. Many of these technical difficulties can be avoided when studying Lévy processes instead of random walks. The results for Lévy processes are completely analogous.

All random elements considered are assumed to be defined on a common probability space $(\Omega, \mathcal{F}, P)$. Denote by $\mathbb{D} = \mathbb{D}([0,1], \mathbb{R}^d)$ the space of càdlàg functions $\mathbf{x} : [0,1] \to \mathbb{R}^d$ equipped with the $J_1$-metric, referred to as $d_0$ as in [3], which makes $\mathbb{D}$ a complete separable linear metric space. In the proofs we will also use the equivalent to $d_0$ incomplete $J_1$-metric, $d$. We denote by $\mathbb{S}_\mathbb{D}$ the "unit sphere" $\{\mathbf{x} \in \mathbb{D} : |\mathbf{x}|_\infty = 1\}$ with $|\mathbf{x}|_\infty = \sup_{t \in [0,1]} |\mathbf{x}_t|$, equipped with the relativized topology of $\mathbb{D}$. Define $\overline{\mathbb{D}}_0 = (0, \infty] \times \mathbb{S}_\mathbb{D}$, where $(0, \infty]$ is equipped with the metric $\rho(x, y) = |1/x - 1/y|$, making it complete and separable. For any element $\mathbf{x} \in \overline{\mathbb{D}}_0$, we write $\mathbf{x} = (x^*, \tilde{\mathbf{x}})$, where $x^* = |\mathbf{x}|_\infty$ and $\tilde{\mathbf{x}} = \mathbf{x}/x^*$. Then $\overline{\mathbb{D}}_0$, equipped with the metric $\max\{\rho(x^*, y^*), d_0(\tilde{\mathbf{x}}, \tilde{\mathbf{y}})\}$, is a complete separable metric space. The topological spaces $\mathbb{D} \setminus \{\mathbf{0}\}$, equipped with the relativized topology of $\mathbb{D}$, and $(0, \infty) \times \mathbb{S}_\mathbb{D}$, equipped with the relativized topology of $\overline{\mathbb{D}}_0$, are homeomorphic; the function $T$ given by $T(\mathbf{x}) = (|\mathbf{x}|_\infty, \mathbf{x}/|\mathbf{x}|_\infty)$ is a homeomorphism. Hence,

$$\mathcal{B}(\overline{\mathbb{D}}_0) \cap ((0, \infty) \times \mathbb{S}_\mathbb{D}) = \mathcal{B}(T(\mathbb{D} \setminus \{\mathbf{0}\})),$$

that is, the Borel sets of $\mathcal{B}(\overline{\mathbb{D}}_0)$ that are of interest to us can be identified with the usual Borel sets on $\mathbb{D}$ (viewed in spherical coordinates) that do not contain the zero function. For notational convenience, we will throughout the paper identify $\mathbb{D}$ with the product space $[0, \infty) \times \mathbb{S}_\mathbb{D}$ so that expressions like $\overline{\mathbb{D}}_0 \setminus \mathbb{D} (= \{\infty\} \times \mathbb{S}_\mathbb{D})$ make sense. We denote by $\mathcal{B}(\overline{\mathbb{D}}_0) \cap \mathbb{D}$ the Borel sets $B \in \mathcal{B}(\overline{\mathbb{D}}_0)$ such that $B \cap (\{\infty\} \times \mathbb{S}_\mathbb{D}) = \varnothing$.

Regular variation on $\mathbb{R}^d$ (for random vectors) is typically formulated in terms of vague convergence on $\mathcal{B}(\overline{\mathbb{R}}^d \setminus \{\mathbf{0}\})$, where $\overline{\mathbb{R}} = \mathbb{R} \cup \{-\infty, \infty\}$. The topology on $\overline{\mathbb{R}}^d \setminus \{\mathbf{0}\}$ is chosen so that $\mathcal{B}(\overline{\mathbb{R}}^d \setminus \{\mathbf{0}\})$ and $\mathcal{B}(\mathbb{R}^d)$ coincide on $\mathbb{R}^d \setminus \{\mathbf{0}\}$. Moreover, $B \in \mathcal{B}(\overline{\mathbb{R}}^d \setminus \{\mathbf{0}\})$ is relatively compact (or bounded) in $\overline{\mathbb{R}}^d \setminus \{\mathbf{0}\}$ if and only if $B \cap \mathbb{R}^d$ is bounded away from $\mathbf{0}$ (i.e., $\mathbf{0} \notin \overline{B \cap \mathbb{R}^d}$) in $\mathbb{R}^d$.

We will see that regular variation on $\mathbb{D}$ is naturally expressed in terms of so-called $\hat{w}$-convergence of boundedly finite measures on $\overline{\mathbb{D}}_0$. A boundedly finite measure assigns finite measure to bounded sets. A sequence of boundedly finite measures $(m_n)_{n \in \mathbb{N}}$ on a complete separable metric space $\mathbb{E}$ converges to $m$ in the $\hat{w}$-topology, $m_n \xrightarrow{\hat{w}} m$, if $m_n(B) \to m(B)$ for every bounded Borel set $B$ with $m(\partial B) = 0$. If the state space $\mathbb{E}$ is locally compact, which $\overline{\mathbb{D}}_0$ is not but $\overline{\mathbb{R}}^d \setminus \{\mathbf{0}\}$ is, then a boundedly finite measure is called a Radon measure, and $\hat{w}$-convergence coincides with vague convergence and we write $m_n \xrightarrow{v} m$. Finally, we notice that if $m_n \xrightarrow{\hat{w}} m$ and $m_n(\mathbb{E}) \to m(\mathbb{E}) < \infty$, then $m_n \xrightarrow{w} m$. For details on $\hat{w}$-, vague and weak convergence, we refer to [6], Appendix 2. See also [14] for details on vague convergence and [24, 25] for relations between vague convergence, point process convergence and regular variation.

We start by defining regular variation of random vectors (see [24, 25, 26]).

DEFINITION 1.1. An $\mathbb{R}^d$-valued random vector $\mathbf{X}$ is said to be regularly varying if there exist a sequence $(a_n)$, $0 < a_n \uparrow \infty$, and a nonnull Radon measure $\mu$ on $\mathcal{B}(\overline{\mathbb{R}}^d \setminus \{\mathbf{0}\})$ with $\mu(\overline{\mathbb{R}}^d \setminus \mathbb{R}^d) = 0$ such that, as $n \to \infty$,

$$n\mathrm{P}(a_n^{-1}\mathbf{X} \in \cdot) \xrightarrow{v} \mu(\cdot) \qquad \text{on } \mathcal{B}(\overline{\mathbb{R}}^d \setminus \{\mathbf{0}\}).$$

We write $\mathbf{X} \in \mathrm{RV}((a_n), \mu, \overline{\mathbb{R}}^d \setminus \{\mathbf{0}\})$.

REMARK 1.1. (i) The limiting measure $\mu$ necessarily obeys a homogeneity property, that is, there exists an $\alpha > 0$ such that $\mu(uB) = u^{-\alpha}\mu(B)$ for every $u > 0$ and $B \in \mathcal{B}(\overline{\mathbb{R}}^d \setminus \{\mathbf{0}\})$. This follows by standard regular variation arguments; see Theorem 1.14 on page 19 in [16]. We then also refer to regular variation of $\mathbf{X}$ with index $\alpha$.

(ii) $\mathbf{X} \in \mathrm{RV}((a_n), \mu, \overline{\mathbb{R}}^d \setminus \{\mathbf{0}\})$ implies that, as $u \to \infty$,

$$\frac{\mathrm{P}(\mathbf{X} \in u\cdot)}{\mathrm{P}(|\mathbf{X}| > u)} \xrightarrow{v} c\mu(\cdot) \qquad \text{on } \mathcal{B}(\overline{\mathbb{R}}^d \setminus \{\mathbf{0}\}),$$

for some $c > 0$. The sequence $(a_n)$ will always be chosen so that $n\mathrm{P}(|\mathbf{X}| > a_n) \to 1$ and, with this choice of $(a_n)$, it follows that $c = 1$ above.

Next we define a heavy-tailed version of large deviation principle.

DEFINITION 1.2. A sequence $(\mathbf{X}^n)$ of $\mathbb{R}^d$-valued random vectors is said to satisfy a heavy-tailed large deviation principle if there exist a sequence $((\gamma_n, \lambda_n))$, $0 < \gamma_n, \lambda_n \uparrow \infty$, and a nonnull Radon measure $\mu$ on $\mathcal{B}(\overline{\mathbb{R}}^d \setminus \{\mathbf{0}\})$ with $\mu(\overline{\mathbb{R}}^d \setminus \mathbb{R}^d) = 0$ such that, as $n \to \infty$,

$$\gamma_n \mathrm{P}(\lambda_n^{-1}\mathbf{X}^n \in \cdot) \xrightarrow{v} \mu(\cdot) \qquad \text{on } \mathcal{B}(\overline{\mathbb{R}}^d \setminus \{\mathbf{0}\}).$$

We write $(\mathbf{X}^n) \in \mathrm{LD}(((\gamma_n, \lambda_n)), \mu, \overline{\mathbb{R}}^d \setminus \{\mathbf{0}\})$.

In this paper we work with functional large deviations for stochastic processes with càdlàg sample paths. The appropriate version of large deviation principle for such processes is as follows.

DEFINITION 1.3. A sequence $(\mathbf{X}^n)$ of stochastic processes with sample paths in $\mathbb{D}$ is said to satisfy a heavy-tailed large deviation principle if there exist a sequence $((\gamma_n, \lambda_n))$, $0 \leq \gamma_n, \lambda_n \uparrow \infty$, and a nonnull boundedly finite measure $m$ on $\mathcal{B}(\overline{\mathbb{D}}_0)$ with $m(\overline{\mathbb{D}}_0 \setminus \mathbb{D}) = 0$ such that, as $n \to \infty$,

$$\gamma_n \mathrm{P}(\lambda_n^{-1} \mathbf{X}^n \in \cdot) \stackrel{\hat{w}}{\to} m(\cdot) \qquad \text{on } \mathcal{B}(\overline{\mathbb{D}}_0).$$

We write $(\mathbf{X}^n) \in \mathrm{LD}(((\gamma_n, \lambda_n)), m, \overline{\mathbb{D}}_0)$.

REMARK 1.2. In [7] a sequence $(\mu_n)$ of measures on a space $\mathbb{E}$ is said to satisfy a large deviation principle if, for all Borel sets $A$,

(1.4)
$$\begin{aligned} -\inf_{x \in A^\circ} I(x) &\leq \liminf_{n \to \infty} c_n \log \mu_n(A) \\ &\leq \limsup_{n \to \infty} c_n \log \mu_n(A) \\ &\leq -\inf_{x \in \overline{A}} I(x), \end{aligned}$$

where $I : \mathbb{E} \to [0, \infty]$ is called a rate function and $c_n \to 0$. The cases of interest are those where $A$ becomes for a large $n$ a rare event with respect to $\mu_n$. Then (1.4) describes the logarithmic behavior of exponentially fast decaying probabilities (as $c_n$ usually goes to zero hyperbolically fast). Nontrivial results require that the underlying distributions have light tails in the sense of a finite moment generating function on a "sizable" part of the parameter space. In this paper we are primarily interested in regularly varying distributions (for which the moment generating function does not exist). If one denotes $\mu_n(A) = \mathrm{P}(\lambda_n^{-1} \mathbf{X}^n \in A)$, then Definition 1.3 can be viewed as describing the nonlogarithmic counterpart of (1.4) for probabilities that decay, typically, hyperbolically fast. However, the precise relation between Definition 1.3 and regular variation is not completely clear at the moment.

The key result we will need is Theorem 1.1 that establishes functional large deviations for certain Markov processes with increments that are not too strongly dependent in the sense that an extreme jump does not trigger further jumps or oscillations of the same magnitude with a nonnegligible probability. We consider strong Markov processes in the sense of Definition 2 in [10], page 56. Let $\mathbf{X} = (\mathbf{X}_t)_{t \in [0, \infty)}$ be a Markov process on $\mathbb{R}^d$ with transition function $P_{u,v}(\mathbf{x}, B)$. For $r \geq 0$, $t \geq 0$ and $B_{\mathbf{x}, r} = \{\mathbf{y} \in \mathbb{R}^d : |\mathbf{y} - \mathbf{x}| < r\}$, define

$$\alpha_r(t) = \sup\{P_{u,v}(\mathbf{x}, B_{\mathbf{x},r}^c) : \mathbf{x} \in \mathbb{R}^d \text{ and } 0 \leq u \leq v \leq t\}.$$

Our weak dependence (in the tails) condition is

(1.5) $$\lim_{n \to \infty} \alpha_{\varepsilon \lambda_n}(n) = 0 \qquad \text{for all } \varepsilon > 0$$

for an appropriate choice of $(\lambda_n)$ with $\lambda_n \uparrow \infty$.

For an $\mathbb{R}^d$-valued stochastic process $\mathbf{X} = (\mathbf{X}_t)_{t \in [0, \infty)}$, we adopt the notation $\mathbf{X}^n = (\mathbf{X}_{nt})_{t \in [0,1]}$ throughout the rest of the paper.

THEOREM 1.1. *Let $\mathbf{X} = (\mathbf{X}_t)_{t \in [0,\infty)}$ be a strong Markov process with sample paths in $\mathbb{D}[0,\infty)$ satisfying (1.5). Suppose there exist a set $T \subset [0,1]$ containing 0 and 1 and all but at most countably many points of $[0,1]$, a sequence $((\gamma_n, \lambda_n))$, $0 < \gamma_n, \lambda_n \uparrow \infty$, and a collection $\{m_t : t \in T\}$ of Radon measures on $\mathcal{B}(\overline{\mathbb{R}}^d \backslash \{\mathbf{0}\})$, with $m_t(\overline{\mathbb{R}}^d \backslash \mathbb{R}^d) = 0$ and with $m_1$ nonnull, such that, as $n \to \infty$,*

$$\gamma_n \mathrm{P}(\lambda_n^{-1} \mathbf{X}_t^n \in \cdot) \xrightarrow{v} m_t(\cdot) \qquad \textit{on } \mathcal{B}(\overline{\mathbb{R}}^d \backslash \{\mathbf{0}\}) \textit{ for every } t \in T,$$

*and, for any $\varepsilon > 0$ and $\eta > 0$, there exists a $\delta > 0$, $\delta$, $1 - \delta \in T$ such that*

(1.6) $\quad m_\delta(B_{\mathbf{0},\varepsilon}^c) - m_0(B_{\mathbf{0},\varepsilon}^c) \leq \eta \quad \textit{and} \quad m_1(B_{\mathbf{0},\varepsilon}^c) - m_{1-\delta}(B_{\mathbf{0},\varepsilon}^c) \leq \eta.$

*Then $(\mathbf{X}^n) \in \mathrm{LD}(((\gamma_n, \lambda_n)), m, \overline{\mathbb{D}}_0)$, where $m$ is uniquely determined by $\{m_t : t \in T\}$. Furthermore, $m(\mathcal{V}_0^c) = 0$, where*

(1.7) $\qquad \mathcal{V}_0 = \{\mathbf{x} \in \mathbb{D} : \mathbf{x} = \mathbf{y} \mathbb{1}_{[v,1]}, v \in [0,1), \mathbf{y} \in \mathbb{R}^d \backslash \{\mathbf{0}\}\}.$

This is a modification of Theorems 13 and 15 in [13] with $(n, a_n)$ replaced by $(\gamma_n, \lambda_n)$. The proof of Theorem 1.1 is essentially identical. Notice that the limiting measure is concentrated on $\mathcal{V}_0$, the set of nonzero right-continuous step functions with exactly one step.

In the next section we specialize to sums of heavy-tailed i.i.d. random vectors and prove a large deviation principle. That result is used in Section 3 to study multivariate ruin probabilities in the heavy-tailed context, and in Section 4 to study long strange segments in the heavy-tailed multivariate context.

**2. Large deviations for a heavy-tailed random walk process.** In this section we show a large deviation principle for a random walk with i.i.d. $\mathbb{R}^d$-valued step sizes $\mathbf{Z}_i$. For a generic element of this sequence, $\mathbf{Z}$, we assume that it is regularly varying: $\mathbf{Z} \in \mathrm{RV}((a_n), \mu, \overline{\mathbb{R}}^d \backslash \{\mathbf{0}\})$. Recall from Remark 1.1 that $\mathbf{Z}$ is then regularly varying for some $\alpha > 0$. We will also write $\mathbf{Z} \in \mathrm{RV}(\alpha, \mu)$.

Consider the random walk process $(\mathbf{S}_n)$ given by

$$\mathbf{S}_0 = \mathbf{0}, \qquad \mathbf{S}_n = \mathbf{Z}_1 + \cdots + \mathbf{Z}_n, \qquad n \geq 1,$$

and write $\mathbf{S}^n = (\mathbf{S}_{[nt]})_{t \in [0,1]}$ for the càdlàg embedding of $(\mathbf{S}_n)$. It is our aim to derive a functional version of the large deviation results of A. V. Nagaev [19, 20], S. V. Nagaev [21] and Cline and Hsing [5], which were mentioned in the Introduction, for the sequence $(\mathbf{S}_n)$.

THEOREM 2.1. *Assume that* $\mathbf{Z} \in \mathrm{RV}(\alpha, \mu)$ *and consider a sequence* $(\lambda_n)$ *such that* $\lambda_n \uparrow \infty$ *and the conditions*

$$\lambda_n^{-1} \mathbf{S}_n \overset{\mathrm{P}}{\to} \mathbf{0}, \qquad\qquad\qquad\qquad\qquad\qquad\qquad \alpha < 2$$
$$\lambda_n^{-1} \mathbf{S}_n \overset{\mathrm{P}}{\to} \mathbf{0}, \qquad \lambda_n/\sqrt{n^{1+\gamma}} \to \infty \quad \textit{for some } \gamma > 0, \qquad \alpha = 2$$
$$\lambda_n^{-1} \mathbf{S}_n \overset{\mathrm{P}}{\to} \mathbf{0}, \qquad \lambda_n/\sqrt{n \log n} \to \infty, \qquad\qquad\qquad \alpha > 2,$$

*hold. Then* $(\mathbf{S}^n) \in \mathrm{LD}(((\gamma_n, \lambda_n)), m, \overline{\mathbb{D}}_0)$, *where* $\gamma_n = [n\mathrm{P}(|\mathbf{Z}| > \lambda_n)]^{-1}$. *Moreover, the measure $m$ satisfies* $m(\mathcal{V}_0^c) = 0$ *and its one-dimensional restrictions satisfy* $m_t = t\mu$ *for* $t \in [0, 1]$.

REMARK 2.1. It follows from the proof of Lemma 12 in [13] that the finite-dimensional restrictions of $m$ satisfy

$$(2.1) \qquad m_{t_1,\ldots,t_k}(A_1 \times \cdots \times A_k) = \sum_{i=1}^{j}(t_i - t_{i-1})\mu(A_i \cap \cdots \cap A_k),$$

$0 = t_0 \leq t_1 \leq \cdots \leq t_k \leq 1$ with $A_1 \times \cdots \times A_k \in \mathcal{B}(\overline{\mathbb{R}}^{dk} \setminus \{\mathbf{0}\})$ and $j = \inf\{i = 1,\ldots,k : \mathbf{0} \notin A_i\}$. Notice that the relation (2.1) is equivalent to the statement

$$(2.2) \qquad\qquad\qquad m = (\mathrm{Leb} \times \mu) \circ T^{-1},$$

where $T : [0,1] \times (\mathbb{R}^d \setminus \{\mathbf{0}\}) \to \mathbb{D}$ is given by $T(t, \mathbf{x}) = \mathbf{x}\mathbb{1}_{[t,1]}(s), 0 \leq s \leq 1$. From here we immediately conclude that the following property of $m$ in spherical coordinates holds. Let

$$\sigma(\cdot) = \mathrm{P}(\{\boldsymbol{\Theta}\mathbb{1}_{[V,1]}(t), t \in [0,1]\} \in \cdot),$$

where $\boldsymbol{\Theta}$ and $V$ are independent, $V$ is uniformly distributed on $(0,1)$ and $\boldsymbol{\Theta}$ is distributed like the spectral measure of $\mathbf{Z}$, that is,

$$\mathrm{P}(\boldsymbol{\Theta} \in \cdot) = \frac{\mu(\{\mathbf{x} : |\mathbf{x}| > 1, \mathbf{x}/|\mathbf{x}| \in \cdot\})}{\mu(\{\mathbf{x} : |\mathbf{x}| > 1\})}.$$

Then for $x > 0$,

$$\frac{m(\{\mathbf{x} \in \mathbb{D} : |\mathbf{x}|_\infty > x, \mathbf{x}/|\mathbf{x}|_\infty \in \cdot\})}{m(\{\mathbf{x} \in \mathbb{D} : |\mathbf{x}|_\infty > 1\})} = x^{-\alpha}\sigma(\cdot).$$

REMARK 2.2. A light-tailed version of functional large deviations for multivariate random walks is Mogulskii's theorem; see [7], page 152.

REMARK 2.3. Under the conditions of the theorem, one can always choose $\lambda_n = cn$ for any positive $c$ if $\alpha \geq 1$ and $\mathrm{E}(\mathbf{Z}) = \mathbf{0}$. If $\alpha \in (0,2)$, an appeal to [22] yields that the conditions (i) $n\mathrm{P}(|\mathbf{Z}| > \lambda_n) \to 0$ and (ii) $n\lambda_n^{-1}\mathrm{E}(\mathbf{Z}\mathbb{1}_{[0,\lambda_n]}(|\mathbf{Z}|)) \to \mathbf{0}$ are necessary and sufficient for $\lambda_n^{-1}\mathbf{S}_n \overset{\mathrm{P}}{\to} \mathbf{0}$. Condition (ii) is satisfied if (i) holds and one of the following conditions holds:

$\alpha \in (0,1)$, or $\alpha = 1$ and $\mathbf{Z}$ is symmetric, or $\alpha \in (1,2)$ and $\mathrm{E}(\mathbf{Z}) = \mathbf{0}$. These conditions are comparable to those in [5] for $\alpha \in (0,2)$. For $\alpha > 2$, the growth condition on $(\lambda_n)$ is slightly more restrictive than in [21], where one can choose $\lambda_n = a\sqrt{n \log n}$ for any $a > \sqrt{\alpha - 2}$, provided $\mathrm{E}(\mathbf{Z}) = \mathbf{0}$. The reason for the more restrictive assumption is that, for our applications, we need convergence on the whole space $\overline{\mathbb{D}}_0$, and this is not guaranteed by the less restrictive assumption.

REMARK 2.4. We mention in passing that the large deviation relation

$$\frac{\mathrm{P}(\lambda_n^{-1}\mathbf{S}_n \in \cdot)}{n\mathrm{P}(|\mathbf{Z}| > \lambda_n)} \xrightarrow{v} \mu(\cdot) \tag{2.3}$$

has a nice interpretation in terms of point process convergence. To see this, rewrite (2.3) as follows:

$$\frac{n}{r_n}\mathrm{P}(a_n^{-1}\mathbf{S}_{r_n} \in \cdot) \xrightarrow{v} \mu(\cdot), \tag{2.4}$$

where, as usual, the sequence $(a_n)$ satisfies $n\mathrm{P}(|\mathbf{Z}| > a_n) \to 1$ and $(r_n)$ is an integer sequence such that $r_n \to \infty$, $r_n/n \to 0$ and $n\mathrm{P}(|\mathbf{Z}| > \lambda_{r_n}) \to 1$. Then (2.4) is equivalent to the following point process convergence result (see [25], Proposition 3.21):

$$N_n = \sum_{i=1}^{[n/r_n]} \delta_{a_n^{-1}(\mathbf{S}_{ir_n} - \mathbf{S}_{(i-1)r_n})} \xrightarrow{d} N, \tag{2.5}$$

where $\delta_\mathbf{x}$ denotes Dirac measure at $\mathbf{x}$, $\xrightarrow{d}$ stands for convergence in distribution in the space $M_p(\overline{\mathbb{R}}^d \setminus \{\mathbf{0}\})$ of point measures on $\overline{\mathbb{R}}^d \setminus \{\mathbf{0}\}$ equipped with the vague topology and $N$ is a Poisson random measure with mean measure $\mu$. Hence, for any $\mu$-continuity set $A$ bounded away from zero, $\mathrm{P}(N_n(A) = 0) \to \mathrm{P}(N(A) = 0) = \exp\{-\mu(A)\}$. In particular, for the componentwise maxima,

$$M_n^{(i)} = \max_{j=1,\ldots,[n/r_n]} (S_{jr_n}^{(i)} - S_{(j-1)r_n}^{(i)}), \qquad i = 1, \ldots, d,$$

and $A = ([0, x_1] \times \cdots \times [0, x_d])^c$, $x_i \geq 0$, $i = 1, \ldots, d$, we have

$$\mathrm{P}(a_n^{-1}M_n^{(1)} \leq x_1, \ldots, a_n^{-1}M_n^{(d)} \leq x_d)$$
$$\to \mathrm{P}(Y_1 \leq x_1, \ldots, Y_d \leq x_d) = \exp\{-\mu(A)\},$$

where $\mathbf{Y}$ is the vector of the component-wise maxima of the points of the limiting Poisson random measure $N$. If $\mu(A) > 0$ for some set $A$ of this type, then a nondegenerate component $Y_i$ of the limiting vector $\mathbf{Y}$ exists and has a Fréchet distribution $P(Y_i \leq x) = \exp\{-cx^{-\alpha}\}$, $x > 0$, for some $c > 0$. The distribution of $\mathbf{Y}$ is one of the multivariate extreme value distributions, see [25], Chapter 5.

Another relation equivalent to (2.4) is given by

$$r_n^{-1} \sum_{i=1}^{n} \delta_{a_n^{-1}(\mathbf{S}_{ir_n} - \mathbf{S}_{(i-1)r_n})} \xrightarrow{\text{P}} \mu,$$

where $\xrightarrow{\text{P}}$ stands for convergence in probability in the space $M_+(\overline{\mathbb{R}}^d \backslash \{\mathbf{0}\})$ of non-negative Radon measures on $\overline{\mathbb{R}}^d \backslash \{\mathbf{0}\}$, see [25], Exercise 3.5.7 and [24]. This result can be interpreted as a "law of large numbers analogue" to the weak convergence result (2.5).

We start with an auxiliary result about the convergence of the one-dimensional distributions. The proof is similar to the proof of the results in [5, 19, 21].

LEMMA 2.1. *Under the conditions of Theorem* 2.1, *for every* $t \geq 0$,

$$\gamma_n \mathrm{P}(\lambda_n^{-1} \mathbf{S}_{[nt]} \in \cdot) \xrightarrow{v} t\mu(\cdot) \qquad \text{on } \mathcal{B}(\overline{\mathbb{R}}^d \backslash \{\mathbf{0}\}).$$

PROOF. We prove the result for $t = 1$, the general case is completely analogous by switching from $\mathbf{S}_n$ to $\mathbf{S}_{[nt]}$. We start with an upper bound for $\gamma_n \mathrm{P}(\lambda_n^{-1} \mathbf{S}_n \in A)$, where $A$ is bounded away from zero and satisfies $\mu(\partial A) = 0$. In what follows we write, for any Borel set $B \subset \overline{\mathbb{R}}^d \backslash \{\mathbf{0}\}$ and $\varepsilon > 0$,

$$B^\varepsilon = \{\mathbf{x} \in \overline{\mathbb{R}}^d \backslash \{\mathbf{0}\} : |\mathbf{y} - \mathbf{x}| \leq \varepsilon, \mathbf{y} \in B\}.$$

Then

$$\mathrm{P}(\lambda_n^{-1} \mathbf{S}_n \in A)$$
$$\leq n\mathrm{P}(\lambda_n^{-1} \mathbf{Z} \in A^\varepsilon) + \mathrm{P}(\lambda_n^{-1} \mathbf{S}_n \in A, \lambda_n^{-1} \mathbf{Z}_i \notin A^\varepsilon \text{ for all } i = 1, \ldots, n)$$
$$\leq n\mathrm{P}(\lambda_n^{-1} \mathbf{Z} \in A^\varepsilon) + \mathrm{P}(\lambda_n^{-1} |\mathbf{S}_n - \mathbf{Z}_i| > \varepsilon \text{ for all } i = 1, \ldots, n)$$
$$= I_1 + I_2.$$

By regular variation of $\mathbf{Z}$, Remark 1.1(i) and since $\mu(\partial A) = 0$, we have

$$\lim_{\varepsilon \downarrow 0} \lim_{n \to \infty} \gamma_n I_1 = \lim_{\varepsilon \downarrow 0} \mu(A^\varepsilon) = \mu(A).$$

Next we show that, for every $\varepsilon > 0$, $\lim_{n \to \infty} \gamma_n I_2 = 0$. We consider the following disjoint partition of $\Omega$ for $\delta > 0$:

$$B_1 = \bigcup_{1 \leq i < j \leq n} \{|\mathbf{Z}_i| > \delta\lambda_n, |\mathbf{Z}_j| > \delta\lambda_n\},$$

$$B_2 = \bigcup_{i=1}^{n} \{|\mathbf{Z}_i| > \delta\lambda_n, |\mathbf{Z}_j| \leq \delta\lambda_n, j \neq i, j = 1, \ldots, n\},$$

$$B_3 = \left\{ \max_{i=1,\ldots,n} |\mathbf{Z}_i| \leq \delta\lambda_n \right\}.$$

Clearly, $\gamma_n \mathrm{P}(B_1) = o(1)$ and

$$\mathrm{P}(\{|\mathbf{S}_n - \mathbf{Z}_i| > \varepsilon\lambda_n \text{ for all } i = 1, \ldots, n\} \cap B_2)$$

$$= \sum_{k=1}^{n} \mathrm{P}(\{|\mathbf{S}_n - \mathbf{Z}_i| > \varepsilon\lambda_n \text{ for all } i = 1, \ldots, n\}$$

$$\cap \{|\mathbf{Z}_k| > \delta\lambda_n, |\mathbf{Z}_j| \leq \delta\lambda_n, j \neq k, j \leq n\})$$

$$\leq \sum_{k=1}^{n} \mathrm{P}(|\mathbf{S}_n - \mathbf{Z}_k| > \varepsilon\lambda_n, |\mathbf{Z}_k| > \delta\lambda_n)$$

$$= \mathrm{P}(|\mathbf{S}_{n-1}| > \varepsilon\lambda_n)[n\mathrm{P}(|\mathbf{Z}| > \delta\lambda_n)]$$

$$= o(\gamma_n^{-1}),$$

where the last equality holds since $\mathbf{Z}$ is regularly varying. As regards $B_3$, we have

$$\mathrm{P}(\{|\mathbf{S}_n - \mathbf{Z}_i| > \varepsilon\lambda_n \text{ for all } i = 1, \ldots, n\} \cap B_3)$$

$$\leq \mathrm{P}\left(|\mathbf{S}_{n-1}| > \varepsilon\lambda_n, \max_{i=1,\ldots,n-1} |\mathbf{Z}_i| \leq \delta\lambda_n\right)$$

$$\leq \sum_{k=1}^{d} \mathrm{P}\left(|S_{n-1}^{(k)}| > \frac{\varepsilon\lambda_n}{d}, \max_{i=1,\ldots,n-1} |Z_i^{(k)}| \leq \delta\lambda_n\right).$$

Therefore, it suffices to show that, for every $k = 1, \ldots, d$ and $\varepsilon > 0$,

$$\mathrm{P}\left(|S_n^{(k)}| > \varepsilon\lambda_n, \max_{i=1,\ldots,n} |Z_i^{(k)}| \leq \delta\lambda_n\right) = o(n\mathrm{P}(|Z^{(k)}| > \lambda_n)).$$

We may assume without loss of generality that $d = 1$ and we adapt the notation correspondingly. Since $\lambda_n^{-1} S_n \xrightarrow{\mathrm{P}} 0$, $n\lambda_n^{-1} \mathrm{E}(Z\mathbb{1}_{[0,\delta\lambda_n]}(|Z|)) \to 0$ for every fixed $\delta > 0$. Hence, for large $n$,

$$\mathrm{P}\left(|S_n| > \varepsilon\lambda_n, \max_{i=1,\ldots,n} |Z_i| \leq \delta\lambda_n\right)$$

$$\leq \mathrm{P}\left(\left|\sum_{i=1}^{n} Z_i \mathbb{1}_{[0,\delta\lambda_n]}(|Z_i|)\right| > \varepsilon\lambda_n\right)$$

$$\leq \mathrm{P}\left(\left|\sum_{i=1}^{n} (Z_i \mathbb{1}_{[0,\delta\lambda_n]}(|Z_i|) - \mathrm{E}(Z\mathbb{1}_{[0,\delta\lambda_n]}(|Z|)))\right| > \frac{\varepsilon\lambda_n}{2}\right).$$

An application of the Fuk–Nagaev inequality (e.g., [22], page 78) yields that the right-hand side is bounded by

$$I_3 = c_1 n \lambda_n^{-p} \operatorname{E}(|Z|^p \mathbb{1}_{[0,\delta\lambda_n]}(|Z|)) + \exp\{-c_2 \lambda_n^2 [n \operatorname{var}(Z \mathbb{1}_{[0,\delta\lambda_n]}(|Z|))]^{-1}\}$$
$$= I_{3,1} + I_{3,2},$$

for any $p \geq 2$, some $c_1, c_2 > 0$. By Karamata's theorem (e.g., [4]), for any $p > \alpha$,

$$\operatorname{E}(|Z|^p \mathbb{1}_{[0,\delta\lambda_n]}(|Z|)) \sim c(\delta\lambda_n)^p \operatorname{P}(|Z| > \delta\lambda_n),$$

as $n \to \infty$. Hence, for $p > \max(2, \alpha)$,

$$\lim_{\delta \downarrow 0} \limsup_{n \to \infty} \frac{\operatorname{E}(|Z|^p \mathbb{1}_{[0,\delta\lambda_n]}(|Z|))}{\lambda_n^p \operatorname{P}(|Z| > \lambda_n)} = c \lim_{\delta \downarrow 0} \limsup_{n \to \infty} \frac{(\delta\lambda_n)^p \operatorname{P}(|Z| > \delta\lambda_n)}{\lambda_n^p \operatorname{P}(|Z| > \lambda_n)}$$
$$= c \lim_{\delta \downarrow 0} \delta^{p-\alpha} = 0.$$

We consider 3 distinct cases to bound $I_{3,2}$:

(i) If $\operatorname{var}(Z) < \infty$, then since $\lambda_n / \sqrt{n \log n} \to \infty$,

$$(2.6) \qquad \limsup_{n \to \infty} \frac{I_{3,2}}{n \operatorname{P}(|Z| > \lambda_n)} = 0.$$

(ii) If $\alpha \in (0, 2)$, by Karamata's theorem,

$$n \lambda_n^{-2} \operatorname{var}(Z \mathbb{1}_{[0,\delta\lambda_n]}(|Z|)) \sim c n \operatorname{P}(|Z| > \lambda_n).$$

Hence, (2.6) holds.

(iii) If $\alpha = 2$ and $\operatorname{var}(Z) = \infty$, then $\operatorname{P}(|Z| > \lambda_n) \lambda_n^2$ and $\operatorname{var}(Z \mathbb{1}_{[0,\delta\lambda_n]}(|Z|))$ are slowly varying functions of $\lambda_n$. Taking into account that $\lambda_n n^{-(1+\gamma)/2} \to \infty$ for some $\gamma > 0$, we conclude that (2.6) holds. We conclude that

$$(2.7) \qquad \limsup_{n \to \infty} \gamma_n \operatorname{P}(\lambda_n^{-1} \mathbf{S}_n \in A) \leq \mu(A^\varepsilon) \to \mu(A) \qquad \text{as } \varepsilon \downarrow 0$$

for any $\mu$-continuity set $A$ bounded away from zero.

To prove the corresponding lower bound, it suffices to consider rectangles $A = [\mathbf{a}, \mathbf{b}) \subset \mathbb{R}^d$ bounded away from zero. These are $\mu$-continuity sets and they determine vague convergence on the Borel $\sigma$-field $\mathcal{B}(\overline{\mathbb{R}}^d \backslash \{\mathbf{0}\})$ by virtue of the fact that $\mu(\overline{\mathbb{R}}^d \backslash \mathbb{R}^d) = 0$. With $\mathbf{a}^{+\varepsilon} = (a_1 + \varepsilon, \ldots, a_d + \varepsilon)'$ and $\mathbf{b}^{-\varepsilon} = (b_1 - \varepsilon, \ldots, b_d - \varepsilon)'$, introduce the set $A^{-\varepsilon} = (\mathbf{a}^{+\varepsilon}, \mathbf{b}^{-\varepsilon}]$, which is a nonempty $\mu$-continuity set for sufficiently small $\varepsilon > 0$. Then

$$\operatorname{P}(\lambda_n^{-1} \mathbf{S}_n \in A) \geq \operatorname{P}(\lambda_n^{-1} \mathbf{S}_n \in A, \lambda_n^{-1} \mathbf{Z}_i \in A^{-\varepsilon} \text{ for some } i \leq n)$$
$$\geq \operatorname{P}(\lambda_n^{-1} \mathbf{Z}_i \in A^{-\varepsilon}, \lambda_n^{-1} |\mathbf{S}_n - \mathbf{Z}_i| < \varepsilon \text{ for some } i \leq n)$$
$$\geq n \operatorname{P}(\lambda_n^{-1} \mathbf{Z} \in A^{-\varepsilon}) \operatorname{P}(\lambda_n^{-1} |\mathbf{S}_{n-1}| < \varepsilon)$$
$$- \frac{n(n-1)}{2} [\operatorname{P}(\lambda_n^{-1} \mathbf{Z} \in A^{-\varepsilon})]^2.$$

Notice that $\mathbf{S}_{n-1}/\lambda_n \xrightarrow{\mathrm{P}} \mathbf{0}$. Hence,

$$
\begin{aligned}
&\liminf_{n\to\infty} \gamma_n \mathrm{P}(\lambda_n^{-1}\mathbf{S}_n \in A) \\
&\qquad \geq \lim_{n\to\infty} \frac{\mathrm{P}(\lambda_n^{-1}\mathbf{Z} \in A^{-\varepsilon})}{\mathrm{P}(|\mathbf{Z}| > \lambda_n)} = \mu(A^{-\varepsilon}) \to \mu(A) \qquad \text{as } \varepsilon \downarrow 0,
\end{aligned}
\tag{2.8}
$$

since $A$ is a $\mu$-continuity set. We conclude from (2.7) and (2.8) that, for every rectangle $A = (\mathbf{a}, \mathbf{b}]$,

$$\lim_{n\to\infty} \gamma_n \mathrm{P}(\lambda_n^{-1}\mathbf{S}_n \in A) = \mu(A).$$

The latter relations determine the vague convergence $\gamma_n \mathrm{P}(\lambda_n^{-1}\mathbf{S}_n \in \cdot) \xrightarrow{v} \mu(\cdot)$. This concludes the proof. $\square$

PROOF OF THEOREM 2.1. It follows immediately from Lemma 2.1 that, for every $t \geq 0$, $\gamma_n \mathrm{P}(\lambda_n^{-1}\mathbf{S}_t^n \in \cdot) \xrightarrow{v} t\mu(\cdot)$. The process $(\mathbf{S}_{[t]})_{t\in[0,\infty)}$ is a strong Markov process satisfying the conditions of Theorem 1.1, which immediately yields that $(\mathbf{S}^n) \in \mathrm{LD}(((\gamma_n, \lambda_n)), m, \overline{\mathbb{D}}_0)$ for some boundedly finite measure $m$ on $\mathcal{B}(\overline{\mathbb{D}}_0)$ satisfying (2.1) and that $m(\mathcal{V}_0^c) = 0$. $\square$

**3. Ruin probabilities for a multivariate random walk with drift.** In this section we are interested in extensions of the notion of ruin probability to an $\mathbb{R}^d$-valued random walk with regularly varying step sizes. We use the same notation as in Section 2, that is, $(\mathbf{Z}_i)$ is an i.i.d. $\mathbb{R}^d$-valued sequence such that $\mathbf{Z} \in \mathrm{RV}(\alpha, \mu)$. Moreover, we assume that $\alpha > 1$. Then $\mathrm{E}(\mathbf{Z})$ is well defined and we assume that $\mathrm{E}(\mathbf{Z}) = \mathbf{0}$. Then we know from Theorem 2.1 that $(\mathbf{S}^n) \in \mathrm{LD}((([n\mathrm{P}(|\mathbf{Z}| > n)]^{-1}, n)), m, \overline{\mathbb{D}}_0)$. We will use this result to derive the asymptotic behavior of the probabilities, as $u \to \infty$,

$$\psi_u(A) = \mathrm{P}(\mathbf{S}_n - \mathbf{c}n \in uA \text{ for some } n \geq 1),$$

$\mathbf{c}$ is a vector and $A$ is a measurable set.

Given $\mathbf{c} \neq \mathbf{0}$, let $\delta > 0$ be such that the set

$$K_{\mathbf{c}}^\delta = \{\mathbf{x} \in \mathbb{R}^d : |\mathbf{x}/|\mathbf{x}| + \mathbf{c}/|\mathbf{c}|| < \delta\}$$

satisfies $\mu((\partial K_{\mathbf{c}}^\delta)\setminus\{\mathbf{0}\}) = 0$. We will take $A \in \mathcal{B}(\mathbb{R}^d \setminus K_{\mathbf{c}}^\delta)$ to avoid sets $A$ that can be hit by simply drifting in the direction $-\mathbf{c}$. Recall from Theorem 2.1 that

$$\gamma_n \mathrm{P}(\mathbf{S}^n \in n\cdot) \xrightarrow{\hat{w}} m(\cdot),$$

where $m$ concentrates on step functions with one step. Using this, we can describe the intuition behind the main result of this section, Theorem 3.1, as follows. Essentially, for large $n$, the random walk process $\mathbf{S}^n$ reaches a set $nA$ for some $t$ by taking one large jump to the set. For the random walk with drift, $\mathbf{S}_{[nt]} - \mathbf{c}[nt]$, the process first drifts in direction $-\mathbf{c}$. Then, at some time $[nv]$, it takes a large jump to a point $-\mathbf{c}[nv] + \mathbf{y}$ and then continues to drift in direction $-\mathbf{c}$. Hence, for $\mathbf{S}_{[nt]} - \mathbf{c}[nt]$, to hit a set $nA$ for some $t$, the jump $\mathbf{y}$ must be of the form $\mathbf{y} = \mathbf{c}[nv] + \mathbf{z} + \mathbf{c}u$, some $\mathbf{z} \in nA$ and $u \geq 0$. That is, $\mathbf{y} \in \mathbf{c}[nv] + \{\mathbf{z} : \mathbf{z} \in \mathbf{c}u + nA, u \geq 0\}$. This explains the appearance of the sets $B_{\mathbf{c}}$ in Theorem 3.1.

Our main result is the following.

THEOREM 3.1. *Assume that* $\mathbf{Z} \in \mathrm{RV}(\alpha, \mu)$ *for some* $\alpha > 1$ *and* $\mathrm{E}(\mathbf{Z}) = \mathbf{0}$. *Then for any set* $A \in \mathcal{B}(\mathbb{R}^d \setminus K_{\mathbf{c}}^{\delta})$ *bounded away from* $\mathbf{0}$,

$$(3.1) \quad \mu^*(A^\circ) \leq \liminf_{u \to \infty} \frac{\psi_u(A)}{u\mathrm{P}(|\mathbf{Z}| > u)} \leq \limsup_{u \to \infty} \frac{\psi_u(A)}{u\mathrm{P}(|\mathbf{Z}| > u)} \leq \mu^*(\overline{A}),$$

*where, for any set* $B \in \mathcal{B}(\mathbb{R}^d \setminus K_{\mathbf{c}}^{\delta})$,

$$\mu^*(B) = \int_0^\infty \mu(\mathbf{c}v + B_{\mathbf{c}})\, dv$$

*and*

$$(3.2) \quad B_{\mathbf{c}} = \{\mathbf{x} + \mathbf{c}t, \mathbf{x} \in B, t \geq 0\}.$$

REMARK 3.1. Notice that neither $\psi_u$ nor $\mu^*$ are additive set functions and, hence, they are not measures. Therefore, (3.1) cannot be stated in terms of vague convergence toward $\mu^*$.

REMARK 3.2. Call a set $A$ $\mathbf{c}$-*increasing* if $\mathbf{x} + \mathbf{c}t \in A$ whenever $\mathbf{x} \in A$ and $t \geq 0$. For such sets, $A_{\mathbf{c}} = A$. If $\mu(\mathbf{c}v + \partial A) = 0$ for almost all $v \geq 0$, then $\mu^*(A^\circ) = \mu^*(\overline{A})$, and Theorem 3.1 gives us

$$\lim_{u \to \infty} \frac{\psi_u(A)}{u\mathrm{P}(|\mathbf{Z}| > u)} = \mu^*(A).$$

An example would be a half space $A = a\mathbf{d} + \{\mathbf{x} : (\mathbf{x}, \mathbf{d}) \geq 0\}$ for some $\mathbf{d}$ with $(\mathbf{d}, \mathbf{c}) > 0$ and $a > 0$. The reason is that, because of the scaling property of the measure $\mu$, it cannot assign a positive mass to any hyperplane unless it contains the origin. Assuming for the ease of notation that $\mathbf{c}$ has positive components, another example is the set $A = \prod_{i=1}^d [x_i, \infty)$ for $\mathbf{x} = (x_1, \ldots, x_d) \in [0, \infty)^d \setminus \{\mathbf{0}\}$.

REMARK 3.3. Notice that the set $B_{\mathbf{c}}$ is universally measurable, and so $\mu^*(B)$ is well defined. Furthermore, it is clear that if $B$ is open, then so is $B_{\mathbf{c}}$. Moreover, if $B$ is closed, then, again, so is $B_{\mathbf{c}}$. To see this, let $\mathbf{y}_n = \mathbf{c}t_n + \mathbf{x}_n \in B_{\mathbf{c}}$ with $t_n \geq 0$ and $\mathbf{x}_n \in B$ for $n = 1, 2, \ldots$. Let $\mathbf{y}_n \to \mathbf{y}$ as $n \to \infty$. If the sequence $(t_n)$ has an accumulation point, it follows from the fact that $B$ is closed that $\mathbf{y} \in B_{\mathbf{c}}$. Therefore, to show that $B_{\mathbf{c}}$ is closed, it is sufficient to show that the sequence $(t_n)$ cannot converge to infinity. Assume, to the contrary, that $t_n \to \infty$. Then

$$\left|\frac{\mathbf{x}_n}{|\mathbf{x}_n|} + \frac{\mathbf{c}}{|\mathbf{c}|}\right| = \left|\frac{\mathbf{y}_n - \mathbf{c}t_n}{|\mathbf{y}_n - \mathbf{c}t_n|} + \frac{\mathbf{c}}{|\mathbf{c}|}\right| = \left|\frac{\mathbf{y}_n/t_n - \mathbf{c}}{|\mathbf{y}_n/t_n - \mathbf{c}|} + \frac{\mathbf{c}}{|\mathbf{c}|}\right| \to \left|\frac{-\mathbf{c}}{|\mathbf{c}|} + \frac{\mathbf{c}}{|\mathbf{c}|}\right| = 0,$$

contradicting the fact that $B \in \mathcal{B}(\mathbb{R}^d \setminus K_{\mathbf{c}}^\delta)$.

REMARK 3.4. In the case $d=1$, relation (3.1) with $A = [1, \infty)$, $\mu(A) > 0$ and $c > 0$ reads as follows:

$$\psi_u(A) = \mathrm{P}\left(\sup_{n \geq 1}(S_n - nc) > u\right) \sim \frac{1}{(\alpha-1)c} u \mathrm{P}(Z > u).$$

This is the classical asymptotic result for the ruin probability in the case of regularly varying $Z_i$'s; see [8] and [9], Chapter 1.

We start the proof with some auxiliary results.

LEMMA 3.1. *For every $A \in \mathcal{B}(\mathbb{R}^d \setminus K_{\mathbf{c}}^\delta)$ bounded away from $\mathbf{0}$,*

$$\lim_{M \to \infty} \limsup_{u \to \infty} \frac{\mathrm{P}(\bigcup_{n > uM} \{\mathbf{S}_n \in n\mathbf{c} + uA\})}{u \mathrm{P}(|\mathbf{Z}| > u)} = 0.$$

PROOF. There exist finitely many points $\mathbf{a}_i$, $i = 1, \ldots, k$, with $(\mathbf{c}, \mathbf{a}_i) > 0$ such that the sets $A_{\mathbf{a}_i} = \{\mathbf{x} \in \mathbb{R}^d : (\mathbf{a}_i, \mathbf{x}) > 1\}$ satisfy $A \subset \bigcup_{i=1}^k A_{\mathbf{a}_i} \cup (\overline{\mathbb{R}}^d \setminus \mathbb{R}^d)$. Hence,

$$\mathrm{P}\left(\bigcup_{n > uM}\{\mathbf{S}_n \in n\mathbf{c} + uA\}\right) \leq \sum_{i=1}^k \sum_{n > uM} \mathrm{P}(\mathbf{S}_n \in n\mathbf{c} + uA_{\mathbf{a}_i}) \tag{3.3}$$
$$= \sum_{i=1}^k \sum_{n > uM} \mathrm{P}((\mathbf{S}_n, \mathbf{a}_i) > n(\mathbf{c}, \mathbf{a}_i) + u).$$

It follows from the uniformity of the large deviation results for one-dimensional centered random walks with regularly varying step sizes (e.g., [5]) that the right-hand side of (3.3) is bounded above by

$$c \sum_{i=1}^k \sum_{n > uM} n \mathrm{P}((\mathbf{Z}, \mathbf{a}_i) > n(\mathbf{c}, \mathbf{a}_i) + u) \leq c_1 \sum_{i=1}^k \int_{uM}^\infty \mathrm{P}((\mathbf{Z}, \mathbf{a}_i) > x(\mathbf{c}, \mathbf{a}_i)) \, dx$$
$$\leq c_2 \sum_{i=1}^k \int_{uM}^\infty \mathrm{P}(|\mathbf{Z}| > x(\mathbf{c}, \mathbf{a}_i)/|\mathbf{a}_i|) \, dx$$
$$\sim c_3 M^{1-\alpha} u \mathrm{P}(|\mathbf{Z}| > u),$$

as $u \to \infty$ ($c, c_1, c_2, c_3 > 0$). In the last step we used Karamata's theorem. This proves the lemma. □

LEMMA 3.2. *If* $(\mathbf{X}^n) \in \mathrm{LD}(((\gamma_n, \lambda_n)), m, \overline{\mathbb{D}}_0)$ *and* $(f_n) \subset \mathbb{D}$ *is a sequence of deterministic functions such that* $f_n \to f$, *then*

$$\gamma_n \mathrm{P}(\lambda_n^{-1} \mathbf{X}^n + f_n - f \in \cdot) \stackrel{\hat{w}}{\to} m(\cdot) \qquad \text{on } \mathcal{B}(\overline{\mathbb{D}}_0).$$

PROOF. Let $A \in \mathcal{B}(\overline{\mathbb{D}}_0)$ be closed and bounded and take $\varepsilon > 0$ small enough such that $A^\varepsilon = \{\mathbf{x} \in \overline{\mathbb{D}}_0 : d_0(\mathbf{x}, A) \leq \varepsilon\}$ is closed and bounded. Since $f_n \to f$, we have $d_0(f_n, f) < \varepsilon$ for $n$ sufficiently large. Hence,

$$\limsup_{n \to \infty} \gamma_n \mathrm{P}(\lambda_n^{-1} \mathbf{X}^n + f_n - f \in A) \leq \limsup_{n \to \infty} \gamma_n \mathrm{P}(\lambda_n^{-1} \mathbf{X}^n \in A^\varepsilon)$$
$$\leq m(A^\varepsilon).$$

Since $A$ is closed, as $\varepsilon \to 0$, $m(A^\varepsilon) \to m(A)$ and the conclusion follows from the Portmanteau theorem. $\square$

PROOF OF THEOREM 3.1. Take $A \in \mathcal{B}(\mathbb{R}^d \setminus K_\mathbf{c}^\delta)$ bounded away from $\mathbf{0}$. We start with an upper bound for $\psi_u(A)$. First notice that, for every $K > 0$,

$$(3.4) \quad \begin{aligned} \psi_u(A) &\leq \mathrm{P}(\mathbf{S}_n - \mathbf{c}n \in u(A \cap \{\mathbf{y} : |\mathbf{y}| \leq K\}) \text{ for some } n \geq 0) \\ &\quad + \mathrm{P}(\mathbf{S}_n - \mathbf{c}n \in u(A \cap \{\mathbf{y} : |\mathbf{y}| > K\}) \text{ for some } n \geq 0) \\ &= \psi_u^{(1)}(A) + \psi_u^{(2)}(A). \end{aligned}$$

Let $\varepsilon > 0$ be small enough so that the set $A^\varepsilon = \{\mathbf{y} \in \mathbb{R}^d : \mathbf{x} \in A, |\mathbf{x} - \mathbf{y}| \leq \varepsilon\}$ is bounded away from the origin and $A^\varepsilon \subset \mathbb{R}^d \setminus K_\mathbf{c}^{\delta/2}$. For all $u \geq \max(2, 2\sqrt{K/\varepsilon})$,

$$\text{if } \mathbf{x} \in u(A \cap \{\mathbf{y} : |\mathbf{y}| \leq K\}) \qquad \text{then } \left| \frac{\mathbf{x}}{u} - \frac{\mathbf{x}}{[u]} \right| \leq \varepsilon,$$

and so $\mathbf{x} \in [u]A^\varepsilon$. Therefore, for $M = 1, 2, \ldots,$

$$(3.5) \quad \begin{aligned} \psi_u^{(1)}(A) &\leq \mathrm{P}(\mathbf{S}_n - \mathbf{c}n \in [u]A^\varepsilon \text{ for some } n \geq 0) \\ &\leq \mathrm{P}\left( \bigcup_{n \leq [u]M} \{\mathbf{S}_n \in (n\mathbf{c} + [u]A^\varepsilon)\} \right) \\ &\quad + \mathrm{P}\left( \bigcup_{n > [u]M} \{\mathbf{S}_n \in (n\mathbf{c} + [u]A^\varepsilon)\} \right) \\ &= \psi_u^{(11)}(A) + \psi_u^{(12)}(A). \end{aligned}$$

We have

$$\psi_u^{(11)}(A) \leq \mathrm{P}((M[u])^{-1}(\mathbf{S}_{[M[u]t]} - \mathbf{c}[M[u]t]) \in M^{-1}A^\varepsilon$$
$$\text{for some rational } t \in [0, 1]).$$

Let $f(t) = \mathbf{c}t$ and for a set $E \in \mathcal{B}(\mathbb{R}^d)$,

$$(3.6) \qquad B_E = \{\mathbf{x} \in \mathbb{D} : \mathbf{x}_t \in M^{-1}E \text{ for some rational } t \in [0, 1]\}.$$

Notice that $B_{A^\varepsilon}$ is bounded away from $\mathbf{0}$ in $\mathbb{D}$ since $A^\varepsilon$ is bounded away from $\mathbf{0}$ in $\mathbb{R}^d$. Hence, also $\overline{B}_{A^\varepsilon}$ is bounded away from $\mathbf{0}$. Since $f(t) = \mathbf{c}t$ and $A^\varepsilon \subset \mathbb{R}^d \setminus K_{\mathbf{c}}^{\delta/2}$, also $f + \overline{B}_{A^\varepsilon}$ is bounded away from $\mathbf{0}$ (i.e., bounded in $\overline{\mathbb{D}}_0$). An application of Theorem 2.1, Lemma 3.2 and the Portmanteau theorem yields

$$\text{(3.7)} \quad \limsup_{u \to \infty} \frac{\psi_u^{(11)}(A)}{Mu\mathrm{P}(|\mathbf{Z}| > Mu)} \leq m(f + \overline{B}_{A^\varepsilon})$$
$$= \int_0^1 \mu(\mathbf{y} : \mathbf{y}\mathbb{1}_{[v,1]} \in f + \overline{B}_{A^\varepsilon})\, dv,$$

where at the last step we used (2.2).

Suppose that, for some $\mathbf{y} \in \mathbb{R}^d \setminus \{\mathbf{0}\}$ and $0 < v < 1$, we have $\mathbf{y}\mathbb{1}_{[v,1]} \in f + \overline{B}_{A^\varepsilon}$. Then there are $\mathbf{x}_n \in f + B_{A^\varepsilon}$ and strictly increasing continuous time changes $h_n : [0,1] \to [0,1]$, $h_n(0) = 0, h_n(1) = 1$ for $n \geq 1$ such that

$$\text{(3.8)} \quad \lim_{n \to \infty} \sup_{0 \leq t \leq 1} |\mathbf{y}\mathbb{1}_{[h_n^{-1}(v),1]}(t) - \mathbf{x}_n(t)| = 0$$

and

$$\lim_{n \to \infty} \sup_{0 \leq t \leq 1} |h_n(t) - t| = 0.$$

Let $0 \leq t_n \leq 1$ and $\mathbf{z}_n \in M^{-1}A^\varepsilon$ be such that $\mathbf{x}_n(t_n) = \mathbf{c}t_n + \mathbf{z}_n, n = 1, 2, \ldots$. It follows from the fact that $A^\varepsilon$ is both bounded away from the origin and $A^\varepsilon \subset \mathbb{R}^d \setminus K_{\mathbf{c}}^\delta$ that the sequence of the norms $|\mathbf{c}t_n + \mathbf{z}_n|, n \geq 1$ is bounded away from zero. We conclude from (3.8) that, for all $n$ large enough, we must have $t_n \geq h_n^{-1}(v)$. If $t_*$ is any accumulation point of the sequence $(t_n)$, it follows that $t_* \geq v$. If $t_{n_k} \to t_*$ as $k \to \infty$, then,

$$|\mathbf{y} - (\mathbf{c}t_* + \mathbf{z}_{n_k})| \leq |\mathbf{y} - (\mathbf{c}t_{n_k} + \mathbf{z}_{n_k})| + |\mathbf{c}||t_{n_k} - t_*| \to 0.$$

Therefore, $\mathbf{y} - \mathbf{c}t_* \in M^{-1}\overline{A^\varepsilon} = M^{-1}A^\varepsilon$, and so

$$\int_0^1 \mu(\mathbf{y} : \mathbf{y}\mathbb{1}_{[v,1]} \in f + \overline{B}_{A^\varepsilon})\, dv$$
$$\leq \int_0^1 \mu(\mathbf{y} : \mathbf{y} \in \mathbf{c}t + M^{-1}A^\varepsilon \text{ for some } t \in [v,1])\, dv$$
$$= M^\alpha \int_0^1 \mu(\mathbf{y} : \mathbf{y} \in \mathbf{c}tM + A^\varepsilon \text{ for some } t \in [v,1])\, dv$$
$$= M^{\alpha-1} \int_0^M \mu(\mathbf{y} : \mathbf{y} \in \mathbf{c}t + A^\varepsilon \text{ for some } t \in [v,M])\, dv.$$

Hence, by (3.7),

$$\limsup_{u \to \infty} \frac{\psi_u^{(11)}(A)}{Mu\mathrm{P}(|\mathbf{Z}| > Mu)} \leq M^{\alpha-1} \int_0^M \mu(\mathbf{y} : \mathbf{y} \in \mathbf{c}t + A^\varepsilon \text{ for some } t \in [v,M])\, dv.$$

Letting $M \to \infty$ and using Lemma 3.1 for $\psi_u^{(12)}(A)$, we conclude that, for all $\varepsilon > 0$,

$$(3.9) \quad \limsup_{u \to \infty} \frac{\psi_u^{(1)}(A)}{u\mathrm{P}(|\mathbf{Z}| > u)} \leq \int_0^\infty \mu(\mathbf{y} : \mathbf{y} \in \mathbf{c}t + A^\varepsilon \text{ for some } t \geq v) \, dv.$$

Fix $v > 0$, let $\varepsilon_n \downarrow 0$, and assume

$$\mathbf{y}_0 \in \bigcap_{n=1}^\infty \{\mathbf{y} : \mathbf{y} \in \mathbf{c}t + A^{\varepsilon_n} \text{ for some } t \geq v\}.$$

Then for every $n \geq 1$, we can write $\mathbf{y}_0 = \mathbf{c}t_n + \mathbf{x}_n$ for some $t_n \geq v$ and $\mathbf{x}_n \in A^{\varepsilon_n}$. The sequence $(t_n)$ must be bounded since $A^{\varepsilon_n} \subset \mathbb{R}^d \setminus K_\mathbf{c}^{\delta/2}$ for all $n$ large enough; see the discussion in Remark 3.3. Let $(n_k)$ be a subsequence such that $t_{n_k} \to t_* \geq v$ as $k \to \infty$. Then $\mathbf{x}_{n_k} \to \mathbf{x}_* \in \overline{A}$ as $k \to \infty$ and, hence,

$$\mathbf{y}_0 = \mathbf{c}t_* + \mathbf{x}_{n_k} + \mathbf{c}(t_{n_k} - t_*) \in \mathbf{c}t_* + \overline{A}.$$

Therefore, letting $\varepsilon \downarrow 0$ in (3.9), we conclude that

$$(3.10) \quad \limsup_{u \to \infty} \frac{\psi_u^{(1)}(A)}{u\mathrm{P}(|\mathbf{Z}| > u)} \leq \int_0^\infty \mu(\mathbf{y} : \mathbf{y} \in \mathbf{c}t + \overline{A} \text{ for some } t \geq v) \, dv$$
$$= \int_0^\infty \mu(\mathbf{c}v + (\overline{A})_\mathbf{c}) \, dv = \mu^*(\overline{A}).$$

Furthermore,

$$\psi_u^{(2)}(A) \leq \mathrm{P}(\mathbf{S}_n - \mathbf{c}n \in u((K_\mathbf{c}^\delta)^c \cap \{\mathbf{y} : |\mathbf{y}| > K\}) \text{ for some } n \geq 0)$$
$$\leq \mathrm{P}(\mathbf{S}_n - \mathbf{c}n \in [u]((K_\mathbf{c}^\delta)^c \cap \{\mathbf{y} : |\mathbf{y}| > K\}) \text{ for some } n \geq 0).$$

The argument leading to (3.9) now gives us

$$\limsup_{u \to \infty} \frac{\psi_u^{(2)}(A)}{u\mathrm{P}(|\mathbf{Z}| > u)}$$
$$\leq \int_0^\infty \mu(\mathbf{y} : \mathbf{y} \in \mathbf{c}t + ((K_\mathbf{c}^\delta)^c \cap \{\mathbf{z} : |\mathbf{z}| \geq K\}) \text{ for some } t \geq v) \, dv.$$

Let $0 < \theta < |\mathbf{c}|\delta/2$. Suppose that there is a number $t > 0$ such that there exists $\mathbf{y} \in \mathbf{c}t + (K_\mathbf{c}^\delta)^c$ with $|\mathbf{y}| \leq \theta t$. Let $\mathbf{z} = \mathbf{y} - \mathbf{c}t$. Then

$$\left| \frac{\mathbf{z}}{|\mathbf{z}|} + \frac{\mathbf{c}}{|\mathbf{c}|} \right| = \left| \frac{\mathbf{y} - \mathbf{c}t}{|\mathbf{y} - \mathbf{c}t|} + \frac{\mathbf{c}}{|\mathbf{c}|} \right| \leq \frac{2|\mathbf{y}|}{t|\mathbf{c}|} \leq \frac{2t\theta}{t|\mathbf{c}|} < \delta$$

by the choice of $\theta$, contradicting the fact that $\mathbf{z} \in (K_\mathbf{c}^\delta)^c$. We conclude that

$$(3.11) \quad \begin{aligned} &\int_0^\infty \mu(\mathbf{y} : \mathbf{y} \in \mathbf{c}t + ((K_\mathbf{c}^\delta)^c \cap \{\mathbf{z} : |\mathbf{z}| \geq K\}) \text{ for some } t \geq v) \, dv \\ &\leq \int_0^\infty \mu(\mathbf{y} : |\mathbf{y}| > \theta v, \\ &\qquad \mathbf{y} \in \mathbf{c}t + ((K_\mathbf{c}^\delta)^c \cap \{\mathbf{z} : |\mathbf{z}| \geq K\}) \text{ for some } t \geq v) \, dv \end{aligned}$$

and the integral is finite. Indeed,

$$\{\mathbf{y} : \mathbf{y} \in \mathbf{c}t + ((K_{\mathbf{c}}^{\delta})^c \cap \{\mathbf{z} : |\mathbf{z}| \geq K\}) \text{ for some } t \geq v\} \subset \{\mathbf{z} : |\mathbf{z}| \geq \delta' K\},$$

with $\delta' = \delta/2$. Hence,

$$\int_0^\infty \mu(\mathbf{y} : |\mathbf{y}| > \theta v, \mathbf{y} \in \mathbf{c}t + ((K_{\mathbf{c}}^{\delta})^c \cap \{\mathbf{z} : |\mathbf{z}| \geq K\}) \text{ for some } t \geq v)\, dv$$

$$\leq \int_0^{\delta' K/\theta} \mu(\mathbf{z} : |\mathbf{z}| > \delta' K)\, dv + \int_{\delta' K/\theta}^\infty (\theta v)^{-\alpha} \mu(\mathbf{z} : |\mathbf{z}| > 1)\, dv$$

$$= (\delta' K)^{1-\alpha} \mu(\mathbf{y} : |\mathbf{y}| > 1) \frac{\alpha}{\theta(\alpha - 1)} \to 0,$$

as $K \to \infty$, which establishes the upper bound in (3.1).

To prove the lower bound in the theorem, notice that, for every $K > 0$ and all $\varepsilon > 0$ small enough, the argument we used to establish (3.5) shows that

$$\psi_u(A) \geq \mathrm{P}(\mathbf{S}_n - \mathbf{c}n \in [u](A_\varepsilon \cap \{\mathbf{y} : |\mathbf{y}| \leq K\}) \text{ for some } n \geq 0)$$

for all $u$ large enough, where $A_\varepsilon = \{\mathbf{x} \in A : \mathbf{y} \in A \text{ for all } \mathbf{y} \text{ with } |\mathbf{y} - \mathbf{x}| < \varepsilon\}$. Denoting $D_{\varepsilon, K} = A_\varepsilon \cap \{\mathbf{y} : |\mathbf{y}| \leq K\}$ and using the notation in (3.6), we conclude by Theorem 2.1, Lemma 3.2 and the Portmanteau theorem that, for every $M = 1, 2, \ldots$,

$$(3.12) \quad \liminf_{u \to \infty} \frac{\psi_u(A)}{Mu\mathrm{P}(|\mathbf{Z}| > Mu)} \geq m(f + B_{D_{\varepsilon,K}}^\circ)$$
$$= \int_0^1 \mu(\mathbf{y} : \mathbf{y}\mathbb{1}_{[v,1]} \in f + B_{D_{\varepsilon,K}}^\circ)\, dv.$$

Again, fix a set $E$ and suppose that, for some $\mathbf{y} \in \mathbb{R}^d \setminus \{\mathbf{0}\}$ and $0 < v < 1$, we have $\mathbf{y} - \mathbf{c}t_* \in M^{-1} E^\circ$ for some $t_* \in [v, 1]$. Let us check that

$$(3.13) \quad \mathbf{y}\mathbb{1}_{[v,1]} \in f + B_E^\circ.$$

To this end, select $\delta > 0$ small enough so that $\{\mathbf{z} : |\mathbf{y} - \mathbf{c}t_* - \mathbf{z}| < \delta\} \subset M^{-1} E^\circ$, and consider any function $\mathbf{x}$ such that

$$(3.14) \quad d(\mathbf{y}\mathbb{1}_{[v,1]}, \mathbf{x}) < \frac{\delta}{3}\left(1 \wedge \frac{1}{|\mathbf{c}|}\right),$$

where $d$ refers to the incomplete Skorohod $J_1$-metric. Let $h$ be a strictly increasing continuous time change, $h : [0, 1] \to [0, 1]$, $h(0) = 0, h(1) = 1$ such that

$$|h(t) - t| < \frac{\delta}{2}\left(1 \wedge \frac{1}{|\mathbf{c}|}\right) \quad \text{and} \quad |\mathbf{y}\mathbb{1}_{[v,1]}(t) - \mathbf{x}(h(t))| < \frac{\delta}{2}\left(1 \wedge \frac{1}{|\mathbf{c}|}\right)$$

for all $0 \leq t \leq 1$. In particular,

$$|\mathbf{y} - \mathbf{x}(h(t_*))| \leq \frac{\delta}{2}\left(1 \wedge \frac{1}{|\mathbf{c}|}\right),$$

so that

$$|(\mathbf{y} - \mathbf{c}t_*) - (\mathbf{x}(h(t_*)) - \mathbf{c}h(t_*))| < \delta.$$

If $h(t_*) = 1$, this already tells us by the choice of $\delta$ that $\mathbf{x} \in f + B_E$. If $h(t_*) < 1$, select a rational $t_0 \in [h(t_*), 1]$ such that

$$|(\mathbf{y} - \mathbf{c}t_*) - (\mathbf{x}(t_0) - \mathbf{c}t_0)| < \delta,$$

implying once again that $\mathbf{x} \in f + B_E$. Therefore, any $\mathbf{x}$ satisfying (3.14) is in $f + B_E$, and so (3.13) holds. We conclude that

$$\int_0^1 \mu(\mathbf{y} : \mathbf{y}\mathbb{1}_{[v,1]} \in f + B^\circ_{D_{\varepsilon,K}})\,dv$$
$$\geq \int_0^1 \mu(\mathbf{y} : \mathbf{y}\mathbb{1}_{[v,1]} \in f + M^{-1}D^\circ_{\varepsilon,K})\,dv$$
$$= M^{\alpha-1}\int_0^M \mu(\mathbf{y} : \mathbf{y} \in \mathbf{c}t + D^\circ_{\varepsilon,K} \text{ for some } t \in [v,M])\,dv.$$

Letting $M \to \infty$, we conclude by (3.12) that

$$\liminf_{u\to\infty} \frac{\psi_u(A)}{u\mathrm{P}(|\mathbf{Z}| > u)} \geq \int_0^\infty \mu(\mathbf{y} : \mathbf{y} \in \mathbf{c}t + D^\circ_{\varepsilon,K} \text{ for some } t \geq v)\,dv.$$

Letting first $K \to \infty$ and then $\varepsilon \to 0$, we conclude that

$$\liminf_{u\to\infty} \frac{\psi_u(A)}{u\mathrm{P}(|\mathbf{Z}| > u)} \geq \int_0^\infty \mu(\mathbf{y} : \mathbf{y} \in \mathbf{c}t + A^\circ \text{ for some } t \geq v)\,dv = \mu^*(A^\circ),$$

establishing the lower bound in (3.1). $\square$

**4. Long strange segments.** In this section we study the notion of long strange segments of $\mathbb{R}^d$-valued random walks with regularly varying steps. Let $(\mathbf{Z}_i)$ be an i.i.d. sequence of random vectors in $\mathbb{R}^d$, and $\mathbf{S}_0 = \mathbf{0}$, $\mathbf{S}_n = \mathbf{Z}_1 + \cdots + \mathbf{Z}_n$, $n \geq 1$.

For a set $A \in \mathcal{B}(\mathbb{R}^d)$ bounded away from $\mathbf{0}$, let

$$R_n(A) = \sup\{k : \mathbf{S}_{i+k} - \mathbf{S}_i \in kA \text{ for some } i \in \{0,\ldots,n-k\}\}.$$

Since we are dealing with the intervals over which the sample mean is "far away" from the true mean, the random variable $R_n(A)$ is often called the length of the long strange segment, or long rare segment. See, for example, [7]. The following theorem describes the large deviations of $R_n(A)$ in the heavy-tailed case. It can be motivated as follows. Suppose first that the set $A$ is increasing (i.e., $t\mathbf{x} \in A$ for all $\mathbf{x} \in A$ and $t \geq 1$). We know from Theorem 2.1 that, for large $n$, $\mathbf{S}^n$ may be approximated by a step function with one step. The long strange segment is therefore due to the large jump. If $R_n(A) > nt$, then the large jump must fall in the set $ntA$, which is essentially the same as saying $\mathbf{S}_n \in ntA$. Hence, for large $n$,

$$\frac{\mathrm{P}(R_n(A) > nt)}{n\mathrm{P}(|\mathbf{Z}| > n)} \approx \frac{\mathrm{P}(\mathbf{S}_n \in ntA)}{n\mathrm{P}(|\mathbf{Z}| > n)} \to \mu(tA).$$

For $A$ nonincreasing, we need to be a bit more careful. To handle this case, we define, for any $A \in \mathcal{B}(\mathbb{R}^d)$ and $0 \leq t < 1$,

(4.1) $$A^*(t) = \bigcup_{t \leq s \leq 1} s\overline{A}, \qquad A^\circ(t) = \bigcup_{t < s \leq 1} sA^\circ.$$

Notice that $A^*(t)$ is a closed set and $A^\circ(t)$ is an open set.

THEOREM 4.1. *Suppose $\mathbf{Z} \in \mathrm{RV}(\alpha, \mu)$ for some $\alpha > 1$ and $\mathrm{E}(\mathbf{Z}) = \mathbf{0}$. Then, for every $t \in (0,1)$ and $A \in \mathcal{B}(\mathbb{R}^d)$ bounded away from $\mathbf{0}$,*

$$\mu(A^\circ(t)) \leq \liminf_{n \to \infty} \frac{\mathrm{P}(n^{-1}R_n(A) > t)}{n\mathrm{P}(|\mathbf{Z}| > n)} \leq \limsup_{n \to \infty} \frac{\mathrm{P}(n^{-1}R_n(A) > t)}{n\mathrm{P}(|\mathbf{Z}| > n)} \leq \mu(A^*(t)).$$

REMARK 4.1. Obviously, if $\mathrm{E}(\mathbf{Z}) = \mathbf{z}$ and $A \in \mathcal{B}(\mathbb{R}^d)$ bounded away from $\mathbf{z}$, then

$$\mu((A - \mathbf{z})^\circ(t)) \leq \liminf_{n \to \infty} \frac{\mathrm{P}(n^{-1}R_n(A) > t)}{n\mathrm{P}(|\mathbf{Z}| > n)} \leq \limsup_{n \to \infty} \frac{\mathrm{P}(n^{-1}R_n(A) > t)}{n\mathrm{P}(|\mathbf{Z}| > n)}$$
$$\leq \mu((A - \mathbf{z})^*(t)).$$

REMARK 4.2. If the set $A$ is *increasing*, then it is easy to check that $A^*(t) = t\overline{A}$ and $A^\circ(t) = tA^\circ$ for all $0 < t < 1$, in which case the scaling property of the measure $\mu$ allows us to state the conclusion of Theorem 4.1 as

$$t^{-\alpha}\mu(A^\circ) \leq \liminf_{n \to \infty} \frac{\mathrm{P}(n^{-1}R_n(A) > t)}{n\mathrm{P}(|\mathbf{Z}| > n)} \leq \limsup_{n \to \infty} \frac{\mathrm{P}(n^{-1}R_n(A) > t)}{n\mathrm{P}(|\mathbf{Z}| > n)} \leq t^{-\alpha}\mu(\overline{A}).$$

For the proof of Theorem 4.1, we need two technical lemmas. For a given set $A \in \mathcal{B}(\mathbb{R}^d)$, let $h_A : \mathbb{D} \to [0,1]$ be given by

$$h_A(\mathbf{x}) = \sup\{t \in [0,1] : \mathbf{x}(s+t) - \mathbf{x}(s) \in tA \text{ for some } s \in [0, 1-t]\}$$

with the convention $\sup \varnothing = 0$. Recall the definition of $\mathcal{V}_0$ from (1.7).

LEMMA 4.1. *Let $A \in \mathcal{B}(\mathbb{R}^d)$ be bounded away from $\mathbf{0}$. If $t \in (0,1)$, then:*

(1) $h_{A^\circ}^{-1}((t,1])$ *is open,*

(2) $\mathcal{V}_0 \cap \overline{h_A^{-1}((t,1])} \subset \mathcal{V}_0 \cap h_{\overline{A}}^{-1}([t,1])$.

PROOF. We first show (1). If $A^\circ = \varnothing$, then $h_{A^\circ}^{-1}((t,1]) = \varnothing$. Therefore, we can assume that $A^\circ \neq \varnothing$. Take $\mathbf{y} \in h_{A^\circ}^{-1}((t,1])$. Then there exists $t^* > t$ and $s \in [0, 1-t^*]$ such that $\mathbf{y}(t^*+s) - \mathbf{y}(s) \in t^*A^\circ$. Since $A^\circ$ is open, there exists $\delta > 0$ such that $\{\mathbf{x} \colon |(\mathbf{y}(t^*+s) - \mathbf{y}(s))/t^* - \mathbf{x}| < \delta\} \subset A^\circ$. Let, once again, $d$ be the incomplete Skorohod metric on the space $\mathbb{D}$, and for a small $\delta' > 0$, let $d(\mathbf{z}, \mathbf{y}) < \delta'$. Let $h$ be a strictly increasing continuous time change, $h \colon [0,1] \to [0,1]$, $h(0) = 0, h(1) = 1$ such that

$$|h(t) - t| < 2\delta' \quad \text{and} \quad |\mathbf{y}(t) - \mathbf{z}(h(t))| < 2\delta' \qquad \text{for all } 0 \leq t \leq 1.$$

Notice that, in particular, $t^* - 4\delta' \leq h(t^*+s) - h(s) \leq t^* + 4\delta'$. Therefore,

$$\left| \frac{\mathbf{z}(h(t^*+s)) - \mathbf{z}(h(s))}{h(t^*+s) - h(s)} - \frac{\mathbf{y}(t^*+s) - \mathbf{y}(s)}{t^*} \right|$$

$$\leq |\mathbf{y}(t^*+s) - \mathbf{y}(s)| \left| \frac{1}{t^*} - \frac{1}{h(t^*+s) - h(s)} \right|$$

$$+ \frac{1}{h(t^*+s) - h(s)} |(\mathbf{z}(h(t^*+s)) - \mathbf{z}(h(s))) - (\mathbf{y}(t^*+s) - \mathbf{y}(s))|$$

$$\leq \frac{4\delta'}{(t^* - 4\delta')} \left( \frac{|\mathbf{y}(t^*+s) - \mathbf{y}(s)|}{t^*} + 1 \right) < \delta$$

if $\delta'$ is small enough. By the choice of $\delta$, this implies that $\mathbf{z}(h(t^*+s)) - \mathbf{z}(h(s)) \in (h(t^*+s) - h(s))A^\circ$, and so

$$h_{A^\circ}(\mathbf{z}) \geq h(t^*+s) - h(s) > t^* - 4\delta' > t$$

if $\delta'$ is small enough. Hence, $\mathbf{z} \in h_{A^\circ}^{-1}((t,1])$, and the latter set is open.

We now show (2). Let $(\mathbf{x}_n)$ be a sequence of elements in $h_A^{-1}((t,1])$ such that $\mathbf{x}_n \to \mathbf{x}$ for some $\mathbf{x} = \mathbf{y}\mathbb{1}_{[v,1]} \in \mathcal{V}_0$. For $n \geq 1$, let $t_n > t$ and $s_n \in [0, 1-t_n]$ be such that

$$\frac{\mathbf{x}_n(s_n + t_n) - \mathbf{x}_n(s_n)}{t_n} \in A.$$

Since $\mathbf{x}_n \to \mathbf{x}$, there exists a sequence $(\lambda_n)$ of strictly increasing continuous mappings of $[0,1]$ onto itself satisfying $\sup_{s \in [0,1]} |\lambda_n(s) - s| \to 0$ and

$$\sup_{s \in [0,1]} |\mathbf{x}_n(s) - \mathbf{x}(\lambda_n(s))| \to 0$$

as $n \to \infty$. In particular, for every $\delta > 0$, there exists $N(\delta)$ such that, for $n > N(\delta)$,

$$\sup_{s \in [0,1]} |\lambda_n(s) - s| < \delta, \qquad \sup_{s \in [0,1]} |\mathbf{x}_n(s) - \mathbf{x}(\lambda_n(s))| < \delta.$$

Take any $\varepsilon, \varepsilon' > 0$. Then, uniformly in $n > N(\delta)$,

$$\left|\frac{\mathbf{x}(\lambda_n(t_n+s_n)) - \mathbf{x}(\lambda_n(s_n))}{\lambda_n(t_n+s_n) - \lambda_n(s_n)} - \frac{\mathbf{x}_n(s_n+t_n) - \mathbf{x}_n(s_n)}{t_n}\right|$$

$$\leq |\mathbf{x}(\lambda_n(t_n+s_n)) - \mathbf{x}(\lambda_n(s_n))|\left|\frac{1}{\lambda_n(t_n+s_n) - \lambda_n(s_n)} - \frac{1}{t_n}\right|$$

$$+ \frac{1}{t_n}|(\mathbf{x}(\lambda_n(t_n+s_n)) - \mathbf{x}(\lambda_n(s_n))) - (\mathbf{x}_n(s_n+t_n) - \mathbf{x}_n(s_n))|$$

$$\leq \frac{2\delta}{t_n}\left(\frac{2|\mathbf{y}|}{(t_n - 2\delta)} + 1\right) < \varepsilon,$$

if $\delta$ is small enough. Therefore,

$$\frac{\mathbf{x}(\lambda_n(t_n+s_n)) - \mathbf{x}(\lambda_n(s_n))}{\lambda_n(t_n+s_n) - \lambda_n(s_n)} \in A^\varepsilon.$$

If $\varepsilon$ is so small that $A^\varepsilon$ is bounded away from $\mathbf{0}$, we conclude that

$$\frac{\mathbf{y}}{\lambda_n(t_n+s_n) - \lambda_n(s_n)} \in A^\varepsilon$$

for all $n$ large enough. Since for $n$ large enough, $\lambda_n(t_n+s_n) - \lambda_n(s_n) \geq t - \varepsilon'$, we conclude that, for all $\varepsilon, \varepsilon' > 0$, $h_{A^\varepsilon}(\mathbf{x}) \geq t - \varepsilon'$. Letting $\varepsilon' \to 0$, we see that, for any $\varepsilon > 0$, $h_{A^\varepsilon}(\mathbf{x}) \geq t$. By letting $\varepsilon \to 0$, we conclude that $\mathbf{x} \in h_A^{-1}([t,1])$. □

LEMMA 4.2. *Let $\delta \in (0,1)$. Then*

(4.2) $$\{n^{-1}R_n(A) > \delta\} \subset \{h_A(n^{-1}\mathbf{S}^n) > \delta\}.$$

*Furthermore, if $\sup_{\mathbf{x} \in A} |\mathbf{x}| < \infty$, then, for every $\varepsilon > 0$ and $1 > \delta' > \delta$,*

(4.3) $$\{n^{-1}R_n(A) > \delta\} \supset \{h_{A_\varepsilon}(n^{-1}\mathbf{S}^n) > \delta'\}$$

*for all $n$ large enough, where $A_\varepsilon = \{\mathbf{x} \in A : \mathbf{y} \in A \text{ for all } \mathbf{y} \text{ with } |\mathbf{y} - \mathbf{x}| < \varepsilon\}$.*

PROOF. Suppose that $n^{-1}R_n(A) = n^{-1}k > \delta$. Then there exist $i \in \{0,\ldots, n-k\}$ such that $\mathbf{S}_{k+i} - \mathbf{S}_i \in kA$. Take $t = n^{-1}k$ and $s = n^{-1}i$. Then

$$n^{-1}(\mathbf{S}_{[n(t+s)]} - \mathbf{S}_{[ns]}) \in tA,$$

that is, $h_A(n^{-1}\mathbf{S}^n) \geq n^{-1}k > \delta$.

In the opposite direction, let $t \in (\delta', 1]$ and $s \in [0, 1-t]$ be such that $n^{-1}(\mathbf{S}_{[n(t+s)]} - \mathbf{S}_{[ns]}) \in tA_\varepsilon$. Then the assumption $\sup_{\mathbf{x} \in A} |\mathbf{x}| < \infty$ implies that

$$\frac{\mathbf{S}_{[n(t+s)]} - \mathbf{S}_{[ns]}}{[n(t+s)] - [ns]} \in \frac{nt}{[n(t+s)] - [ns]} A_\varepsilon \subset A$$

for all $n$ large enough, and so

$$R_n(A) \geq [n(t+s)] - [ns] > nt - 1 > n\delta' - 1 > n\delta$$

for all $n$ large enough. $\square$

PROOF OF THEOREM 4.1. Take $t \in (0,1)$, and $A \in \mathcal{B}(\mathbb{R}^d)$ bounded away from $\mathbf{0}$. By Theorem 2.1, $(\mathbf{S}^n) \in \mathrm{LD}(((\gamma_n, \lambda_n)), m, \overline{\mathbb{D}}_0)$ with $\lambda_n = n$ and $\gamma_n = [n\mathrm{P}(|\mathbf{Z}| > n)]^{-1}$. Since $m(\mathcal{V}_0^c) = 0$,

$$m \circ h_{\overline{A}}^{-1}([t,1]) = \mathrm{Leb} \times \mu(\{(v, \mathbf{y}) \in [0,1] \times \mathbb{R}^d : h_{\overline{A}}(\mathbf{y}\mathbb{1}_{[v,1]}) \in [t,1]\})$$

$$= \mathrm{Leb} \times \mu(\{(v, \mathbf{y}) \in [0,1] \times \mathbb{R}^d : \mathbf{y} \in s\overline{A} \text{ for some } t \leq s \leq 1\})$$

$$= \mu(A^*(t)).$$

Therefore, by Lemma 4.2, the Portmanteau theorem and Lemma 4.1(2), we have

$$\limsup_{n \to \infty} \frac{\mathrm{P}(n^{-1}R_n(A) > t)}{n\mathrm{P}(|\mathbf{Z}| > n)} \leq \limsup_{n \to \infty} \frac{\mathrm{P}(h_A(n^{-1}\mathbf{S}^n) > t)}{n\mathrm{P}(|\mathbf{Z}| > n)}$$

$$\leq \limsup_{n \to \infty} \frac{\mathrm{P}(n^{-1}\mathbf{S}^n \in \overline{h_A^{-1}((t,1])})}{n\mathrm{P}(|\mathbf{Z}| > n)}$$

$$\leq m(\overline{h_A^{-1}((t,1])})$$

$$\leq m(h_{\overline{A}}^{-1}([t,1]))$$

$$= \mu(A^*(t)),$$

thus, establishing the upper bound in the theorem.

For the lower bound, suppose first that $\sup_{\mathbf{x} \in A} |\mathbf{x}| \leq C$ for some $C < \infty$. Then by Lemma 4.2, the Portmanteau theorem and Lemma 4.1(1), we have, for every $\varepsilon > 0$ and $t' \in (t,1]$,

$$\liminf_{n \to \infty} \frac{\mathrm{P}(n^{-1}R_n(A) > t)}{n\mathrm{P}(|\mathbf{Z}| > n)}$$

$$\geq \liminf_{n \to \infty} \frac{\mathrm{P}(h_{A_\varepsilon}(n^{-1}\mathbf{S}^n) > t')}{n\mathrm{P}(|\mathbf{Z}| > n)}$$

$$\geq \liminf_{n \to \infty} \frac{\mathrm{P}(h_{(A_\varepsilon)^\circ}(n^{-1}\mathbf{S}^n) > t')}{n\mathrm{P}(|\mathbf{Z}| > n)}$$

$$\geq m(h_{(A_\varepsilon)^\circ}^{-1}(t',1])$$

$$= \mathrm{Leb} \times \mu(\{(v, \mathbf{y}) \in [0,1] \times \mathbb{R}^d : h_{(A_\varepsilon)^\circ}(\mathbf{y}\mathbb{1}_{[v,1]}) > t'\})$$

$$= \mathrm{Leb} \times \mu(\{(v, \mathbf{y}) \in [0,1] \times \mathbb{R}^d : \mathbf{y} \in s(A_\varepsilon)^\circ \text{ for some } t < s \leq 1\})$$

$$= \mu\left(\bigcup_{t'<s\leq 1} s(A_\varepsilon)^\circ\right).$$

Letting first $t' \downarrow t$ and then $\varepsilon \downarrow 0$, we conclude that

$$\liminf_{n\to\infty} \frac{\mathrm{P}(n^{-1}R_n(A) > t)}{n\mathrm{P}(|\mathbf{Z}| > n)} \geq \mu\left(\bigcup_{t<s\leq 1} sA^\circ\right),$$

hence, establishing the lower bound in the theorem for sets $A$ bounded in $\mathbb{R}^d$. In the general case, let, for $C > 0$, $A_{(C)} = \{\mathbf{x} \in A : |\mathbf{x}| \leq C\}$. Then by what we already know,

$$\liminf_{n\to\infty} \frac{\mathrm{P}(n^{-1}R_n(A) > t)}{n\mathrm{P}(|\mathbf{Z}| > n)} \geq \liminf_{n\to\infty} \frac{\mathrm{P}(n^{-1}R_n(A_{(C)}) > t)}{n\mathrm{P}(|\mathbf{Z}| > n)} \geq \mu\left(\bigcup_{t<s\leq 1} sA_{(C)}^\circ\right),$$

and by letting $C \uparrow \infty$, we obtain

$$\liminf_{n\to\infty} \frac{\mathrm{P}(n^{-1}R_n(A) > t)}{n\mathrm{P}(|\mathbf{Z}| > n)} \geq \mu\left(\bigcup_{t<s\leq 1} sA^\circ\right) = \mu(A^\circ(t)),$$

as required. $\square$

In conclusion we derive the distributional limit of the length $R_n(A)$ of long strange segments under a different, nonlarge-deviation, scaling. Let $a_n$ be an increasing sequence such that

(4.4) $$n\mathrm{P}(|\mathbf{Z}| > a_n) \to 1 \quad \text{as } n \to \infty.$$

Notice that $a_n$ is regularly varying with index $1/\alpha$.

THEOREM 4.2. *Suppose $\mathbf{Z} \in \mathrm{RV}(\alpha, \mu)$ for some $\alpha > 1$ and $\mathrm{E}(\mathbf{Z}) = \mathbf{0}$. Then for every $A \in \mathcal{B}(\mathbb{R}^d)$ bounded away from $\mathbf{0}$ and every $x > 0$,*

$$\exp\left\{-x^{-\alpha}\mu\left(\bigcup_{s\geq 1} s\overline{A}\right)\right\} \leq \liminf_{n\to\infty} \mathrm{P}(a_n^{-1}R_n(A) \leq x)$$

$$\leq \limsup_{n\to\infty} \mathrm{P}(a_n^{-1}R_n(A) \leq x)$$

$$\leq \exp\left\{-x^{-\alpha}\mu\left(\bigcup_{s\geq 1} sA^\circ\right)\right\}.$$

In particular, if $\mu(\bigcup_{s\geq 1} sA^\circ) = \mu(\bigcup_{s\geq 1} s\overline{A}) := v$, then

(4.5) $$a_n^{-1}R_n(A) \xrightarrow{d} v^{1/\alpha}W,$$

where $W$ is a standard Fréchet random variable with distribution $\mathrm{P}(W \leq w) = e^{-w^{-\alpha}}, w > 0$.

REMARK 4.3. For the asymptotic behavior of $R_n(A)$ in the light tailed case, see [7], Theorem 3.2.1. In the heavy-tailed case, one-dimensional versions of (4.5) are well known, and not only in the i.i.d. case. See [17] and [23].

REMARK 4.4. If the set $A$ is increasing (see Remark 4.2), then the result of the theorem can be stated in the form

$$\exp\{-x^{-\alpha}\mu(\overline{A})\} \leq \liminf_{n\to\infty} P(a_n^{-1}R_n(A) \leq x)$$
$$\leq \limsup_{n\to\infty} P(a_n^{-1}R_n(A) \leq x)$$
$$\leq \exp\{-x^{-\alpha}\mu(A^\circ)\},$$

and the weak convergence in (4.5) holds whenever $A$ is a $\mu$-continuity set, in which case $v = \mu(A)$.

PROOF OF THEOREM 4.2. Observe that, for every $n \geq k$ and $t > 0$ by independence,

$$(4.6) \qquad P(R_n(A) \leq t) \leq (P(R_k(A) \leq t))^{[n/k]}.$$

Selecting $t = xa_n$ and $k = [Ma_n]$ for $M > x$, we obtain, by (4.6),

$$P(a_n^{-1}R_n(A) \leq x) \leq (P(R_{[Ma_n]}(A) \leq a_n x))^{[n/[Ma_n]]}$$
$$\leq \left[1 - P\left(\frac{1}{[Ma_n]}R_{[Ma_n]}(A) > \frac{a_n x}{[Ma_n]}\right)\right]^{(n/Ma_n)-1}.$$

Next, we use the lower bound in Theorem 4.1, the scaling property of the measure $\mu$, the definition of $a_n$ and regular variation to see that, for every $0 < \varepsilon < \min(1, M/x - 1)$, we have, for all $n$ large enough,

$$P(a_n^{-1}R_n(A) \leq x)$$
$$\leq \left[1 - P\left(\frac{1}{[Ma_n]}R_{[Ma_n]}(A) > (1+\varepsilon)\frac{x}{M}\right)\right]^{(n/Ma_n)-1}$$
$$\leq \left[1 - (1-\varepsilon)[Ma_n]P(|\mathbf{Z}| > [Ma_n])\mu\left(\bigcup_{(1+\varepsilon)x/M < s \leq 1} sA^\circ\right)\right]^{n/(Ma_n)-1}$$
$$\sim \left[1 - (1-\varepsilon)[Ma_n]\frac{M^{-\alpha}}{n}M^\alpha(1+\varepsilon)^{-\alpha}\mu\left(\bigcup_{x < s \leq M/(1+\varepsilon)} sA^\circ\right)\right]^{n/(Ma_n)}$$
$$\sim \left[1 - \frac{1-\varepsilon}{(1+\varepsilon)^\alpha}\frac{Ma_n}{n}\mu\left(\bigcup_{x < s \leq M/(1+\varepsilon)} sA^\circ\right)\right]^{n/(Ma_n)}$$

$$\to \exp\left\{-\frac{1-\varepsilon}{(1+\varepsilon)^\alpha}\mu\left(\bigcup_{x<s\leq M/(1+\varepsilon)} sA^\circ\right)\right\}$$

as $n \to \infty$. Letting $\varepsilon \downarrow 0$ and $M \uparrow \infty$, we conclude by the scaling property of $\mu$ that

$$\limsup_{n\to\infty} P(a_n^{-1} R_n(A) \leq x) \leq \exp\left\{-\mu\left(\bigcup_{x<s<\infty} sA^\circ\right)\right\}$$
$$= \exp\left\{-x^{-\alpha}\mu\left(\bigcup_{s\geq 1} sA^\circ\right)\right\},$$

thus, obtaining the upper bound of the theorem.

We now switch to proving the lower bound of the theorem. To this end, notice that, for every $n \geq k$ and $t > 0$,

(4.7)
$$\{R_n(A) > t\} \subset \Bigg\{\text{for some } j = 1,\ldots, \left[\frac{n}{k}\right]+1, \frac{\mathbf{Z}_{i_1+1}+\cdots+\mathbf{Z}_{i_1+i_2}}{i_2} \in A$$
$$\text{for some } (j-1)k \leq i_1 < jk, i_2 > t \text{ and } i_1 + i_2 \leq jk,$$
$$\text{or for some } j = 1,\ldots, \left[\frac{n}{k}\right]+1, \text{ the point } jk \text{ belongs to an}$$
$$\text{interval } (i_1+1, i_1+i_2) \text{ with } i_2 > t \text{ and } \frac{\mathbf{Z}_{i_1+1}+\cdots+\mathbf{Z}_{i_1+i_2}}{i_2} \in A\Bigg\}.$$

We implicitly assume that we have an infinite sequence $(\mathbf{Z}_k)$ and so having a subscript $k > n$ does not cause a problem. As before, we select $t = xa_n$ and $k = [Ma_n]$, this time for some $M > C > x$. The role of the extra parameter $C$ is seen below. We obtain, by (4.7),

$$P(x < a_n^{-1} R_n(A) \leq C)$$
$$\leq P\left(R_{[Ma_n]}^{(i)}(A) > a_n x \text{ for some } i = 1,\ldots, \left[\frac{n}{[Ma_n]}\right]+1\right)$$
$$+ P\left(R_{2[Ca_n]}^{(i)}(A) > a_n x \text{ for some } i = 1,\ldots, \left[\frac{n}{[Ma_n]}\right]+1\right),$$

where $R_k^{(i)}(A)$, $i = 1, 2, \ldots$, are i.i.d. copies of $R_k(A)$. Repeating the argument in the first part of the proof, and using this time the upper bound in Theorem 4.1, we see that

$$\lim_{n\to\infty} P\left(R_{[Ma_n]}^{(i)}(A) > a_n x \text{ for some } i = 1,\ldots, \left[\frac{n}{[Ma_n]}\right]+1\right)$$
$$= 1 - \lim_{n\to\infty}\left[1 - (Ma_n)\frac{M^{-\alpha}}{n}M^\alpha\mu\left(\bigcup_{x\leq s\leq M} s\overline{A}\right)\right]^{n/(Ma_n)}$$

$$= 1 - \exp\left\{-\mu\left(\bigcup_{x \leq s \leq M} s\overline{A}\right)\right\}$$

and

$$\lim_{n \to \infty} P\left(R^{(i)}_{2[Ca_n]}(A) > a_n x \text{ for some } i = 1, \ldots, \left[\frac{n}{[Ma_n]}\right] + 1\right)$$

$$= 1 - \lim_{n \to \infty}\left[1 - (2Ca_n)\frac{(2C)^{-\alpha}}{n}(2C)^\alpha \mu\left(\bigcup_{x \leq s \leq 2C} s\overline{A}\right)\right]^{n/(Ma_n)}$$

$$= 1 - \exp\left\{-\frac{2C}{M}\mu\left(\bigcup_{x \leq s \leq M} s\overline{A}\right)\right\}.$$

Letting $M \to \infty$, we obtain

$$\limsup_{n \to \infty} P(x < a_n^{-1} R_n(A) \leq C) \leq 1 - \exp\left\{-\mu\left(\bigcup_{s \geq x} s\overline{A}\right)\right\}$$

for every $C > x$. Letting now $C \to \infty$, we obtain the required lower bound in the theorem once we show that

(4.8) $$\lim_{C \to \infty} \limsup_{n \to \infty} P(a_n^{-1} R_n(A) > C) = 0.$$

Let $\rho = \inf_{\mathbf{x} \in A} |\mathbf{x}| > 0$, and observe that, for every $t > 0$,

$$\{R_n(A) > t\} \subset \bigcup_{j=1}^d \{R_{n,j}([-\rho/\sqrt{d}, \rho/\sqrt{d}]^c) > t\},$$

where $R_{n,j}(\cdot)$ is the long strange segment corresponding to the $j$th marginal random walk $(S_n^{(j)})$, $j = 1, \ldots, d$. Therefore, by the one-dimensional results (see, e.g., [17]),

$$\limsup_{n \to \infty} P(a_n^{-1} R_n(A) > C) \leq \lim_{n \to \infty} \sum_{j=1}^d P(a_n^{-1} R_{n,j}([-\rho/\sqrt{d}, \rho/\sqrt{d}]^c) > C)$$

$$= \sum_{j=1}^d (1 - \exp\{-K_j C^{-\alpha}\}),$$

where $K_1, \ldots, K_d$ are finite nonnegative numbers, from which (4.8) follows immediately. $\square$

**Acknowledgments.** The final version of this paper was written at the Mittag–Leffler Institute, Djursholm, in October 2004. We thank the organizers of the semester on "Queueing Theory and Teletraffic Theory." We would also like to thank the anonymous referees for their comments that led to an improvement of presentation of the material.

H. Hult
School of Operations Research
  and Industrial Engineering
Cornell University
414A Rhodes Hall
Ithaca, New York 14853
USA
e-mail: hult@orie.cornell.edu
url: www.orie.cornell.edu/~hult/

T. Mikosch
Department of Applied Mathematics
  and Statistics
University of Copenhagen
DK-2100 Copenhagen
Denmark
e-mail: mikosch@math.ku.dk
url: www.math.ku.dk/~mikosch/

F. Lindskog
Department of Mathematics
KTH
SE-100 44 Stockholm
Sweden
e-mail: lindskog@math.kth.se
url: www.math.kth.se/~lindskog/

G. Samorodnitsky
School of Operations Research
  and Industrial Engineering
Cornell University
220 Rhodes Hall
Ithaca, New York 14853
USA
e-mail: gennady@orie.cornell.edu
url: www.orie.cornell.edu/~gennady/